\documentclass[aap,seceqn,nameyear,MSNbibl,noautosecdot,dvips]{arximspdf}
\usepackage{graphicx}

\doi{10.1214/09-AAP657}
\volume{20}
\issue{5}
\pubyear{2010}
\firstpage{1567}
\lastpage{1606}

\makeatletter

\newtheorem{theorem}{Theorem}

\newproclaim{definition}{Definition}[section]
\newproclaim{remark}[definition]{Remark}
\newtheorem{proposition}[definition]{Proposition}
\newtheorem{corollary}[definition]{Corollary}

\makeatother

\begin{document}
\begin{frontmatter}

\title{The diversity of a distributed genome in bacterial
populations\thanksref{T2}}
\runtitle{The distributed genome in bacterial populations}

\thankstext{T2}{Supported by the BMBF through the Freiburg
Initiative for Systems Biology Grant 0313921.}

\begin{aug}
\author[A]{\fnms{F.} \snm{Baumdicker}},
\author[B]{\fnms{W. R.} \snm{Hess}} and
\author[A]{\fnms{P.} \snm{Pfaffelhuber}\corref{}\ead[label=e3]{p.p@stochastik.uni-freiburg.de}}
\runauthor{F. Baumdicker, W. R. Hess and P. Pfaffelhuber}
\affiliation{Albert-Ludwigs University}
\address[A]{F. Baumdicker\\
P. Pfaffelhuber\\
Fakult\"{a}t f\"{u}r Mathematik und Physik\\
Albert-Ludwigs University\\
Eckerstra\ss e 1\\
D-79104 Freiburg\\
Germany\\
\printead{e3}}
\address[B]{W. R. Hess\\
Fakult\"{a}t f\"{u}r Biologie\\
Albert-Ludwigs University\\
Sch\"{a}nzlestr. 1\\
D-79104 Freiburg\\
Germany}
\end{aug}

\received{\smonth{7} \syear{2009}}
\revised{\smonth{10} \syear{2009}}

%
\begin{abstract}
The distributed genome hypothesis states that the set of genes in a
population of bacteria is distributed over all individuals that
belong to the specific taxon. It implies that certain genes can be
gained and lost from generation to generation. We use the random
genealogy given by a Kingman coalescent in order to superimpose
events of gene gain and loss along ancestral lines. Gene gains occur
at a constant rate along ancestral lines. We assume that gained genes
have never been present in the population before. Gene losses occur
at a rate proportional to the number of genes present along the
ancestral line. In this \textit{infinitely many genes model} we derive
moments for several statistics within a sample: the average number
of genes per individual, the average number of genes differing
between individuals, the number of incongruent pairs of genes, the
total number of different genes in the sample and the gene frequency
spectrum. We demonstrate that the model gives a reasonable fit with
gene frequency data from marine cyanobacteria.
\end{abstract}

%

\end{frontmatter}

\section{\texorpdfstring{Introduction.}{Introduction}}
Population genetics is dealing with biological diversity of
species. Concepts developed in this area include models for genetic
drift, mutation, selection, recombination and population structure.
These models are applied frequently to eukaryotic species to analyze
their evolutionary history. For prokaryotes these concepts are applied
less frequently and \citet{MaynardSmith1995} even asked, ``\textit{Do
bacteria have population genetics}?''

Usually (by the \textit{biological species concept}), a species is a
reproductively isolated set of individuals. This definition can hardly
be applied to prokaryotes, that is, bacteria and archea. In microbiology,
researchers have developed other approaches, mostly defining a species
via genomic similarity. This similarity is either based on
hybridization of DNA, or on DNA sequences of specific molecules such
as ribosomal RNA or suitable housekeeping genes, known to be highly
conserved (and which are identified by genome sequencing or by a
technique called \textit{multilocus sequence typing}
[\citet{MaidenEtAl1998}]).
The definition of bacterial species is
complicated by the fact that the similarity of bacteria depends on the
considered genomic region, which can be explained by transfer of
genetic material between these bacteria
[\citet{DykhuizenGreen1991}]. Even more extreme, individuals from the
same species carry different genes. For example, one quarter of the
genome of a pathogenic variant of \textit{E. coli} was found to be
absent in a laboratory strain [\citet{Perna2001}]. Such findings
lead to
several new hypotheses: the \textit{core genome hypothesis} argues that
the set of genes common to all bacteria of a species is responsible
for maintaining species-specific phenotypic properties [e.g.,
\citet{Riley2009}]. The \textit{distributed genome hypothesis}
predicts that no single individual comprises the full set of genes of
the bacterial population [e.g., \citet{Ehrlich2005}].

The distributed genome hypothesis is similar to the idea of bacterial
\textit{pangen\-omes}. Taking the different gene content of individuals of
a population into account, the pangenome consists of all different
genes carried by all individuals. The pangenome can be split into the
\textit{core genome}, that is, the set of genes carried by every individual
of the population, and the \textit{dispensable} (also: \textit{auxiliary}
or \textit{flexible} or \textit{contingency}) \textit{genome}
[\citet{pmid16185861}]. The pangenome was first analyzed for pathogenic
strains of \textit{Streptococcus agalactiae} [\citet
{pmid16172379}]. It
was shown that around 80\% of a single genome (i.e., the genome of a
single individual) forms the core genome. However, each fully
sequenced genome carries genes which do not occur in other strains,
suggesting that the core genome represents only a small fraction of
the pangenome. The situation is even more extreme in
\textit{Prochlorococcus} and \textit{Synechococcus}, which are marine
cyanobacteria, where the core genome consists of around 60\% of the
genes found in a single genome [\citet{pmid18159947},
\citet{Dufresne2008GenomeBiol18507822}]. In contrast,
\citet{pmid16185861} show that a set of four genomes of
\textit{Bacillus
anthracis} contain all genes found in the complete sample of 8
individuals, showing that the core genome is the biggest part of the
pangenome of this species. Recently, the pangenome of all bacteria was
considered, using a dataset of 573 completely sequenced genomes
showing that only 250 genes (which are 8\% of a bacterial genome on
average) were common to almost all bacterial species
[\citet{Lapierre2009}, \citet{Bentley2009}].

The bacterial \textit{supragenome} makes the split into the core and
dispensable genome more precise: each gene present in a population (or
in a sample) has a frequency (for core genes this is 100\%) such that
the pangenome gives rise to a gene frequency spectrum. The first
analysis of \citet{pmid17550610} on a sample of 13 genomes of
\textit{Haemophilus influenzae} shows that the largest class (19\% of
the pangenome) in the dispensable genome are genes only present in a
single genome. In addition, every pair of genomes differs by around
300 genes on average. Similar findings were obtained for
\textit{Streptococcus pneumoniae} [\citet{pmid17675389}].

The pan- and supragenome suggest that genes can be gained and lost
along lineages of bacteria, leading to diversity of genomes. It is
well known that genes can be gained in bacteria by three different
mechanisms: (i) The uptake of genetic material from the environment is
referred to as \textit{transformation}. (ii) Bacteria can be infected by
lysogenic phages which provide additional genetic material that can be
built in the bacterial genome. This process is known as
\textit{transduction}. (iii) A~direct link between two bacterial cells
of the same species leads to exchange of genetic material, known as
\textit{conjugation}. These three mechanisms are usually referred to as
\textit{horizontal gene flow}. Events of gene loss occur by mutations
resulting in pseudogenization or deletion of genes.


The aim of the present paper is to model the bacterial pangenome. We
focus on two different aspects: the genealogical relationships between
individuals and the mutational mechanism. Using the diffusion limit of
a standard neutral model (with finite offspring variance) leads to a
random genealogy, usually referred to as the Kingman coalescent
[\citet{Kingman1982}, \citet{Wakeley2008}]. Gene gain and
loss is the basis of
our mutational model, as introduced by
\citet{Huson2004Bioinformatics15044248} in the phylogenetics
literature. Here, new genes are taken up from the environment at
constant rate along ancestral lines. We assume that all genes taken up
are different. In addition, present genes can be lost at constant
rate. In analogy to standard population genetic models we refer to
this as the \textit{infinitely many genes model}.

\section{\texorpdfstring{Model.}{Model}}\label{S:2}
The dynamics of our model consist of two parts. Reproduction follows
the (diffusion limit of a) neutral Wright--Fisher model (or some other
exchangeable population genetic model with finite offspring
variance). The mutation model we use is borrowed from the
phylogenetics literature
[\citet{Huson2004Bioinformatics15044248}]
and describes gene gains and losses along ancestral lines. After
introducing the model in Sections \ref{S:21} and \ref{S:22}, we
discuss connections to other mutation models in Section \ref{S:23}.

\subsection{\texorpdfstring{Reproduction dynamics.}{Reproduction dynamics}}
\label{S:21}
We will use the neutral Wright--Fisher model:
a~panmictic population of size $N$ reproduces neutrally and clonally,
that is, asexually. In this model, individuals in generation $t+1$ choose
a unique parent from generation $t$ purely at random and independent
of other individuals. It is well known that the genealogy of a sample
of size $n$ taken from the Wright--Fisher model converges for
$N\to\infty$ after a time rescaling by $N$ to the Kingman coalescent
started with $n$ lines [e.g., \citet{Durrett2008},
\citet{Wakeley2008}]. In
this process, starting with $n$ lines:
\begin{itemize}
\item if there are $k$ lines left, draw an exponential time with rate
${k\choose2}$ (which equals the number of pairs in the $k$ lines),
which is the time to the next coalescence event;
\item at the next coalescence event pick two lines at random from the
$k$ lines and merge these into one line.
\end{itemize}
If there is one line left, the sample has found its most recent common
ancestor which we trace back into the past for an infinite amount of
time.
\begin{definition}[(Kingman coalescent)]
We denote the random tree resulting from the above mechanism---the
Kingman coalescent---by $\mathcal T$. We consider $\mathcal T$ as a
partially ordered metric space with order relation $\preceq$ and
metric $d_{\mathcal T}$ where the distance of two points in
$\mathcal T$ is given by the sum of the times to their most recent
common ancestor. We make the convention that $s\preceq t$ for
$s,t\in\mathcal T$ if $s$ is an ancestor of $t$.
\end{definition}

We note that our starting point, the Wright--Fisher model, can be
replaced by other models. For continuous, overlapping generations, the
Moran model is the most canonical choice. Generally, every
exchangeable model with genealogy---under a suitable time
rescaling---given by the Kingman coalescent leads to the same results
as those
obtained in the present paper; the genealogy of a sequence of
exchangeable models converges to the Kingman coalescent if and only if
the sequence of offspring distributions of a single individual has
bounded finite variance and fulfills a condition regarding their third
moments [\citet{MoehleSagitov2001}].\looseness=1

\subsection{\texorpdfstring{Mutation dynamics.}{Mutation dynamics}}
\label{S:22}
We model individuals whose genomes consist of sets of genes. Every
individual has a set of genes $\mathcal G_c$, $g_c:=|\mathcal G_c|$
which are absolutely necessary to survive and hence are conserved,
that is, must be passed from ancestor to offspring. The genes $\mathcal
G_c$ constitute the \textit{core genome}. In addition, we model an
infinite gene pool by a set of genes $I=[0,1]$ with $\mathcal G_c\cap
I = \varnothing$. The genome of individual $i$ in the sample, $1\leq
i\leq N$, contains genes $\mathcal G_i\subseteq I$ which are not
necessary for the individuals to survive. This set of genes is called
the \textit{dispensable genome} of individual $i$.

During the lifetime of every individual or at every reproduction
event, mutations may happen. In our mutation model, the
\textit{infinitely many genes model}, we assume the following two
mechanisms (in a Wright--Fisher population of size $N$) which changes
the dispensable genome from parent to offspring:
\begin{itemize}
\item\textit{gene gain}: before reproduction of individual $i$, there
is a probability $\mu$ that a new gene $u\in I$ is taken up from the
environment. We assume that gene $u$ has never been present in any
genome of the population;
\item\textit{gene loss}: every gene of the dispensable genome
$u\in\mathcal G_i$ of individual $i$ is lost with probability $\nu$
before reproduction of individual $i$.
\end{itemize}
We take an extreme point of view here in that the core genes are
absolutely necessary for an individual to survive and the genes in the
dispensable genome evolve completely neutral. Using this view,
ancestry is not affected by mutations, that is, all gene gains and losses
seen in the population are assumed to be neutral.

Using the same time-scaling as for the genealogies, we assume that
$\mu=\mu_N$ and $\nu=\nu_N$ are such that $\theta= \lim_{N\to
\infty}
2N\mu_N$ and $\rho= \lim_{N\to\infty} 2N\nu_N$. After this rescaling
of the parameters, new genes are gained at rate $\frac\theta2$ and
present genes are lost at rate $\frac\rho2$.
\begin{definition}[(Tree-indexed Markov chain for gene gain and loss)]
\label{def:21}
Let $\mathcal T$ be the Kingman coalescent. We either assume that
$\mathcal T$ is rooted at the most recent common ancestor of the
sample or that $\mathcal T$ has a single infinite line.
Given~$\mathcal T$, we define a Markov chain $\Gamma_{\mathcal T} =
(\mathcal G_t)_{t\in\mathcal T}$, indexed by $\mathcal T$, with
state space $\mathcal N_f(I)$, the space of counting measures on
$I$. [The Markov property for the tree-indexed Markov chain
$\Gamma_{\mathcal T}$ states that for all $t\in\mathcal T$,
$(\mathcal G_s)_{t\preceq s}$ depends on $(\mathcal G_s)_{s\preceq
t}$ only through~$\mathcal G_t$.] Denoting by $\lambda_I$ the
Lebesgue measure on $I$, $\Gamma_{\mathcal T}$ makes
transitions\looseness=1
%
%
\begin{eqnarray}
\label{eq:41}
&&\mbox{from }g \mbox{ to } g + \delta_u \mbox{ at rate }
\frac\theta2
\lambda_I(du),\nonumber\\[-8pt]\\[-8pt]
&&\mbox{from }g \mbox{ to } g - \delta_u \mbox{ at rate }
\frac\rho2 g(du)\nonumber
\end{eqnarray}
along $\mathcal T$. Taking into account that the tree $\mathcal T$
has $n$ leaves, one for each individual of the sample, we denote
these leaves by $1,\ldots,n\in\mathcal T$. In this setting, $\mathcal
G_1,\ldots,\mathcal G_n$ describe the genes present in individuals
$1,\ldots,n$.
\end{definition}

An illustration of the tree-indexed Markov
%
%
%
\begin{figure}[b]

\includegraphics{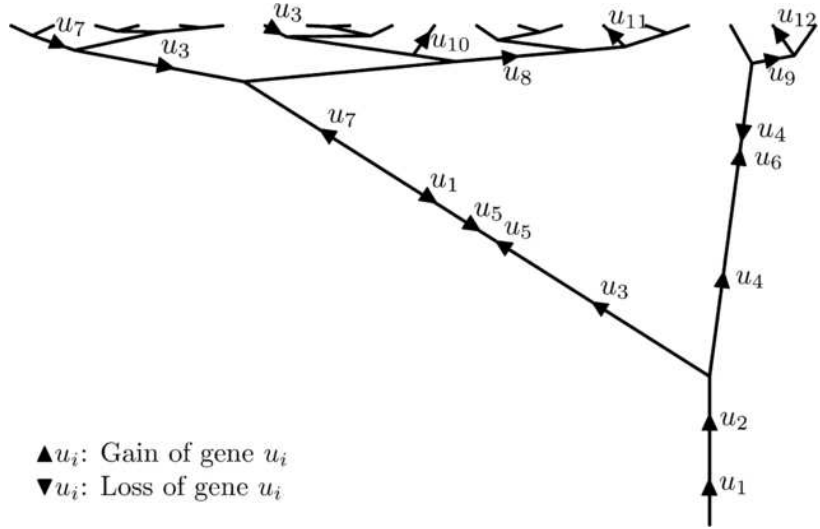}

\caption{An illustration of the infinitely many genes model along a
Kingman coalescent. If a gene is gained along a line
(indicated by the $\blacktriangle$-sign) it can be lost again
(indicated by a $\blacktriangledown$-sign). An individual of the
sample (i.e., a leaf of the coalescent tree) carries the set of
genes which were gained along its ancestral lines and did not get
lost again. Here are some examples: the gene $u_2$ is present in all
individuals, $u_3$ is only present in 10 individuals of the left
branch due to two gene losses. The genes $u_4$ and $u_5$ were lost
in all ancestral lines and do not occur in any individual. The gene
$u_7$ is missing in the 2 individuals on the left-hand side due to a gene
loss and in the right branch as the gene gain was in the left
branch.}
\label{fig:illu}
\end{figure}
chain is shown in Figure \ref{fig:illu}.
\begin{remark}[(Notation)]
\label{rem:finMeas}
Note that all gained genes are almost surely different, so $\mathcal
G_t$ does not have double points, that is, $\mathcal
G_t(\{u\})\in\{0,1\}$, for all\vadjust{\goodbreak} $u\in I$ and $t\in\mathcal T$, almost
surely. We will use the following notation, equating counting
measures without double points with their support: consider
$g\in\mathcal N_f(I)$ without double points. There is $m\in\mathbb
N$ and $u_1,\ldots,u_m$ with $g = \sum_{i=1}^m \delta_{u_i}$. We will
refer to $u_1,\ldots,u_m$ as the points in $g$ and also write
$g=\{u_1,\ldots,u_m\}$. Moreover, we define\looseness=1
\begin{eqnarray*}
|g| &:=& \int1 \,dg,\qquad g_s \cap g_t := g_s\wedge g_t,\\
g_s\setminus g_t &:=& (g_s - g_t)^+
\end{eqnarray*}
for $g,g_s,g_t\in\mathcal
N_f(I)$.
\end{remark}

Our aim is to describe patterns of the dispensable genomes $\mathcal G
:= (\mathcal G_1, \ldots, \mathcal G_n)$ or whole genomes $(\mathcal
G_1\cup\mathcal G_c, \ldots, \mathcal G_n\cup\mathcal G_c )$. These
results can then be compared to genomic data of a sample of bacteria
which gives the genes (or gene families) carried by individuals in the
sample.

\subsection{\texorpdfstring{Comparison to other mutation
models.}{Comparison to other mutation models}}
\label{S:23}
In mathematical population genetics, there are several standard
mutation models, for example, the \textit{infinitely many alleles
model} and
the \textit{infinitely many sites model} [see, e.g., \citet{Durrett2008}
or \citet{Ewens2004}]. In the former, every mutation (along some random
tree) leads to a new type, also called a new allele. It is assumed
that mutated alleles have never been present in the population
before. The latter is a refinement of the former: the allele of an
individual in the sample is modeled as an infinite stretch on
DNA. Every mutation is assumed to change a single site on this genome
(hence leading to a new allele), and it is assumed that every mutation
hits a site that has never been hit before. The last assumption is
relaxed in the \textit{finite sites model} where sites can be hit
several times changing their state between several possibilities.

The \textit{infinitely many genes model} as described above is
conceptually different from these standard models: the infinitely many
sites model (along some random tree $\mathcal T$) can be described,
when a genome is given as the linear set $I$ as above using that
events, occurring at rate $\frac\theta2 \lambda(du)$ along the
tree, changes the state from the ancestral to a derived state at
position $u$ in the genome. However, loss events do not have a
correspondence in the infinitely many sites model. In the finite sites
model, a site can change from the ancestral state to a derived state
and back; however, in the infinitely many genes model, once a gene has
changed from the ancestral state (not present) to the derived state
(present) and back (not present) no further changes of the state are
possible.

Although all mutation models are conceptually different, the infinite
sites model can be seen as the infinitely many genes model for
$\rho=0$. To understand this, consider a random tree with gene gain
events (and no losses due to $\rho=0$). Reinterpreting these gene
gains as point mutations along a chromosome, each\vadjust{\goodbreak} mutation hitting a
new site, leads to the infinitely many sites model. However, there are
still differences between the infinitely many genes model for $\rho=0$
and the infinitely many sites model. On the one hand, $\rho=0$ implies
that genes cannot get lost which leads to an infinite genome for all
individuals. On the other hand, most interesting quantities
concentrate on mutations segregating (i.e., showing both the ancestral
and the mutated---or derived---state) in the sample. In summary, the
infinitely many genes model for $\rho=0$ and the infinitely many sites
model are the same with respect to properties of segregating
genes/sites. However, we will see in our results below (Theorems
\ref{T2}, \ref{T4} and \ref{T5}) that the infinitely many genes model
is not continuous in $\rho=0$ for certain aspects of segregating
mutations.

\section{\texorpdfstring{Results.}{Results}}
In the rest of the paper, we fix a sample of size $n\in\mathbb N$ and
$\theta, \rho>0$. We describe expectations and variances of several
quantities of interest. If we want to stress the dependence on the
model parameters we will use subscripts, for example, $\mathbb
E_{\theta,\rho}[\cdot]$, in order to make this clear. Since the core genome
is conserved for all individuals of the population, we focus on the
dispensable genome first. We provide results for the average number of
genes in the sample (Section \ref{S:31}, Theorem~\ref{T1}), the
average number of pairwise differences (Section \ref{S:32}, Theorem
\ref{T2}), incongruent pairs of genes (Section \ref{S:33}, Theorem
\ref{T3}), the size of the dispensable genome of the sample (Section
\ref{S:34}, Theorem \ref{T4}) and the gene frequency spectrum (Section
\ref{S:35}, Theorem \ref{T5}). We then extend these results to the
complete pangenome, that is, the union of the dispensable and core genome
(Section \ref{S:36}) and describe the application of our model to a
dataset from Prochlorococcus, a marine cyanobacterium (Section
\ref{S:37}). Finally, we discuss biologically realistic extensions of
our model (Section \ref{S:out}).

\subsection{\texorpdfstring{Average number of genes.}{Average number of genes}}
\label{S:31}
The simplest statistics in the infinitely many genes model is based on
counting the number of genes for all individuals in the sample. The
\textit{average number of genes} (\textit{in the dispensable genome})
is given by
%
%
\begin{equation}\label{eq:Hbar}
A := \frac1n \sum_{i=1}^n |\mathcal G_i|.
\end{equation}
Our first results provide the first and second moment for $A$.
\begin{theorem}[(Average number of genes)]\label{T1}
For $A$ as above
\begin{eqnarray*}
\mathbb E[A] & = & \frac{\theta}{\rho},\\
\mathbb V[A] & = & \frac{1}{n} \frac{\theta}{1+\rho} +
\frac{\theta}{\rho(1+\rho)}.
\end{eqnarray*}
\end{theorem}
\begin{remark}
Note that the result for $\mathbb E[A]$ is robust against changes
of the reproduction mechanism or nonequilibrium
situations. Consider any model of reproduction which has not gone
extinct by time $t$. As long as the mutation mechanism is
independent of reproduction, picking an individual at time $t$ from
the population gives a single ancestral line along which genes
accumulate by the same distribution. In particular, the results for
$\mathbb E[A]$ remain unaltered under population size changes or
population subdivision.

The fact that $\mathbb V[A]$ does not converge to 0 as
$n\to\infty$ does not come as a big surprise: the sets $\mathcal
G_1,\ldots,\mathcal G_n$ are dependent through the joint
genealogy---given through the Kingman coalescent---relating the sample.
\end{remark}

\subsection{\texorpdfstring{Average number of pairwise
differences.}{Average number of pairwise differences}}
\label{S:32}
The \textit{average number of pairwise differences} is given by
%
%
\begin{equation}\label{eq:Dbar}
D := {\frac1{n(n-1)} \sum_{1\leq i \neq j\leq n}} |\mathcal G_i
\setminus\mathcal G_j|.
\end{equation}
\begin{theorem}[(Average number of pairwise differences)]\label{T2}
For $D$ as above,
\begin{eqnarray*}
\mathbb E[D] & = &\frac{\theta}{1+\rho},\\
\mathbb V[D] & = &\theta\biggl( \frac{(3 + 14 \rho+ 23 \rho^2 + 16
\rho^3 + 4 \rho^4 + 4 \theta+ 2 \rho\theta)}
{(1+\rho)^2(2+\rho)(3+\rho)(1+2\rho)(3+2\rho)} \\
&&\hspace*{11.46pt}{} + \frac{6 + 19 \rho+ 19 \rho^2 + 12 \rho^3 + 4
\rho^4
+ 8 \theta+ 4 \rho\theta}{(1 + \rho) (2 + \rho) (3 + \rho)
(1 + 2 \rho) (3 + 2 \rho)}\frac1n \\
&&\hspace*{11.46pt}{} + \frac{3 + 11
\rho+ 12 \rho^2 + 4 \rho^3 + 10 \theta+ 9 \rho\theta+ 2
\rho^2 \theta}{(1 + \rho) (2 + \rho) (3 + \rho) (1 + 2 \rho)
(3 + 2 \rho)} \frac{2}{n(n-1)} \biggr).
\end{eqnarray*}
\end{theorem}
\begin{remark}
The quantity $D$ is only based on genes segregating in the
sample. Hence, as explained in Section \ref{S:23}, the infinitely
many genes model for $\rho=0$ is equivalent to the infinitely many
sites model with respect to $D$. As the theorem shows, the expected
number of differences between individuals $i$ and $j$ is $\mathbb
E[|\mathcal G_i\setminus\mathcal G_j| + |\mathcal
G_j\setminus\mathcal G_i|] = 2\frac\theta{1+\rho}$. Hence,\vspace*{1pt}
$\mathbb E_\rho[D]$ is not continuous in $\rho=0$ since the comparable
quantity in the infinite sites model, the average number of
segregating sites in a sample of size two, is $\theta$. The reason
is that for small $\rho$, every individual carries a lot of genes,
all of which can get lost at the small rate $\rho$. These loss
events lead to differences between individuals as well as events of
gene gain. A similar argument shows that the variance is not
continuous [see, e.g., \citet{Wakeley2008}, (4.15)] for the variance
in the infinite sites model.

Note that $\mathbb V[D]$ does not converge to 0 as
$n\to\infty$. Again---the reason is that the differences $(\mathcal
G_i \setminus\mathcal G_j)_{1\leq i\neq j\leq n}$ are dependent
through the underlying common genealogy.
\end{remark}

\subsection{\texorpdfstring{Incongruent pairs of genes.}{Incongruent
pairs of genes}}
\label{S:33}
Assume the following situation: for a pair of genes there are four
individuals in which all four possible states of\vadjust{\goodbreak} presence/absence of
the two genes are observed. This means that the following situation is
found:

\mbox{}

\begin{center}
\begin{tabular}{@{}ccc@{}}
\hline
& gene 1 & gene 2 \\
\hline
Individual 1 & present & present \\
Individual 2 & present & absent \\
Individual 3 & absent & present \\
Individual 4 & absent & absent \\
\hline
\end{tabular}
\end{center}

\mbox{}

\begin{figure}[b]

\includegraphics{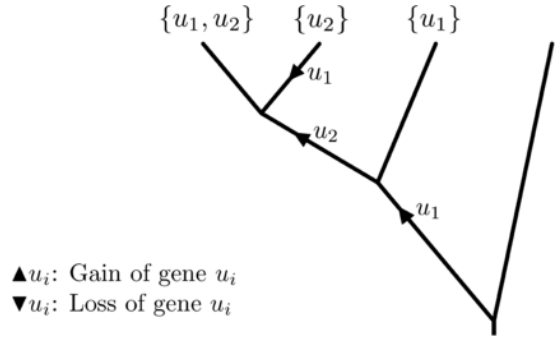}

\caption{If genes can be lost ($\rho>0$) it is
possible that all four possible configurations of presence/absence
of a pair of genes are seen in a sample of four individuals. We call
such pairs incongruent. The figure shows one possible history of the
two genes which leads to incongruence.}\label{fig:four}
\end{figure}

If genes cannot be lost ($\rho=0$) this situation cannot
occur in our model. The reason is that gene 1 would indicate that
individuals 1 and 2 have a common ancestor before 1 and 3 have
(otherwise individual 3 would also carry gene 1), while gene 2
indicates that individuals 1 and 3 have a common ancestor before 1 and
2 have. This is the reason why we call pairs of genes for which the
above situation appears \textit{incongruent}. If $\rho>0$, incongruent
pairs can arise by gene loss; see Figure~\ref{fig:four} for an
example. We will now state how many incongruent pairs we can expect to
see in our sample.

%
%

The \textit{average number of incongruent pairs of genes} (\textit{in four
genomes}) is given by
\[
P := \frac{1}{n(n-1)(n-2)(n-3)} \sum_{i,j,k,l=1}^n D_{ij,kl}\cdot
D_{ik,jl},
\]
where
%
%
\begin{equation}
\label{eq:Aijkl}
D_{ij,kl} := |(\mathcal G_i\cap\mathcal G_j) \setminus(\mathcal
G_k\cup\mathcal G_l)|,\qquad 1\leq i,j,k,l\leq n.
\end{equation}
\begin{theorem}[(Incongruent pairs of genes)]\label{T3}
For $P$ as above,
\begin{eqnarray*}
\mathbb{E}[P] &=& \frac{\theta^2\rho}{4} \biggl({18+117\frac\rho
2+203\frac{\rho^2}{4}+105\frac{\rho^3}{8}}\biggr)\\
&&{}\times\biggl(\biggl(1+\frac\rho
2\biggr)^2\biggl(1+2\frac\rho2\biggr)\biggl(1+4\frac\rho2\biggr)\\
&&\hspace*{18.2pt}{}\times\biggl(3+4\frac\rho
2\biggr)\biggl(3+5\frac\rho2\biggr)\biggl(6+5\frac\rho2\biggr)
\biggl(6+7\frac\rho2\biggr)\biggr)^{-1}.
\end{eqnarray*}
\end{theorem}
\begin{remark}
In the proof of Theorem \ref{T3}, we have to consider all
possible genealogies of four individuals. In order to obtain the
variance of $P$, one would have to take into account all possible
genealogies relating eight individuals.

Note that
%
%
\begin{equation}\label{eq:EAikl}
\mathbb E[D_{ij,kl}] = \frac16 \frac{\theta}{2}\frac{4\cdot
3}{(3+\rho)(2+\rho)} = \frac{\theta}{(3+\rho)(2+\rho)}
\end{equation}
by the gene frequency spectrum (Theorem \ref{T5}). Using this
results it can be shown that $\mathbb{COV}[D_{ij,kl},D_{ik,jl}]<0$
in all cases.

For $\rho=0$ we find that either $D_{ij,kl}=0$ or
$D_{ik,jl}=0$ implying that \mbox{$P=0$} (hence $\mathbb E_{\rho=0}[P]=0$),
almost surely. The theorem\vspace*{-1pt} implies that
$\mathbb{E}_\rho[D_{ij,kl}\cdot D_{ik,jl}]\stackrel{\rho\to
0}{\hbox
to 1cm{\rightarrowfill}}0$,
that is, $\mathbb E_\rho[P]$ is continuous\vspace*{1pt} in $\rho=0$. This
observation is not obvious since small $\rho$ implies that all
individuals carry many genes. Using this fact, one could argue that
the chance to observe a pair of genes giving rise to $D_{ij,kl}\cdot
D_{ik,jl}$ grows with decreasing $\rho$. However, although the
number of genes grows for small $\rho$, pairs of genes giving rise
to $D_{ij,kl}\cdot D_{ik,jl}$ are most probably created by two gene
gains and one gene loss, as shown in Figure \ref{fig:four}, for
small $\rho$.

As we will discuss in Section \ref{S:out}, the
possibility of horizontal gene transfer (by bacterial conjugation)
would be a biologically realistic extension of our model. Under such
a mechanism, new genes are not only taken from the environment, but
also from other individuals of the population. As a result, the
genealogical tree would be different for different genes. Hence, the
order of coalescence can be different, and so there is an increased
number of incongruent pairs of genes. Hence, the theorem is valuable
for determining the possibility of horizontal gene transfer in real
populations.
\end{remark}

\subsection{\texorpdfstring{Size of the dispensable genome.}{Size of
the dispensable genome}}
\label{S:34}
Now we come to properties involving the whole sample (in contrast to
pairs and quartets of individuals in the last theorems). The simplest
statistics involving all individuals of the sample is the total number
of genes, that is, the size of the dispensable genome.

The \textit{size of the dispensable genome} is given by
%
%
\begin{equation}\label{eq:G}
G:= \Biggl| \bigcup_{i=1}^n \mathcal G_i \Biggr|.
\end{equation}

We need the following definition in order to give the variance of the
number of genes in the dispensable genome in our next theorem.
\begin{definition}[(The functions $h_k$ and $g_{\underline k}$)]
\label{def:gh}
For $k\geq0$ we define
%
%
\begin{equation}
\label{eq:defgh0}
h_k := 2\sum_{i=0}^{k-1}\frac{1}{\rho+i}.
\end{equation}
Moreover, for $\underline k = (k_1, k_2, k_3)$ with $k_1, k_2,
k_3\geq0$, we set
\begin{eqnarray*}
\underline k'_1 &=& (k_1-1, k_2, k_3),\qquad \underline k'_2 = (k_1, k_2-1,
k_3), \\
\underline k'_3 &=& (k_1, k_2, k_3-1),\qquad
\underline k'_4 = (k_1+1, k_2-1, k_3-1), \\
\underline k'_5 &=& (k_1-1, k_2+1, k_3),\qquad
\underline k'_6 = (k_1-1, k_2, k_3+1),
\\
\lambda_1 &=& \pmatrix{k_1\cr2},\qquad \lambda_2 = \pmatrix{k_2\cr2}
+ k_1k_2 + \frac\rho2 k_2,\\
\lambda_3 &=& \pmatrix{k_3\cr2}
+ k_1k_3 + \frac\rho2 k_3,\qquad
\lambda_4 = k_2k_3,\qquad \lambda_5 = \lambda_6 = \frac
\rho2k_1
\end{eqnarray*}
and $\overline\lambda= \sum_{i=1}^6\lambda_i$. We define
recursively
%
%
\begin{equation}\label{eq:defgh1}
g_{\underline k} = \cases{ \dfrac2 \rho h_{k_1+k_2}, &\quad if
$k_1+k_3=1$,\cr
\dfrac2 \rho h_{k_1+k_3}, &\quad if $k_1+k_2=1$,}
\end{equation}
and
%
%
\begin{eqnarray}\label{eq:defgh2}
g_{\underline k} &=& (k_1+k_2)(k_1+k_3)\frac{2}{\overline\lambda{}^2}
\\
&&{} + \sum_{i=1}^6 \frac{\lambda_i}{\overline
\lambda} \biggl( \frac1
{\overline\lambda} \bigl((k_1+k_2)h_{k_1'+k_3'} +
(k_1+k_3)h_{k_1'+k_2'} \bigr) + g_{\underline k_i'} \biggr)
\end{eqnarray}
in all other cases.
\end{definition}
\begin{theorem}[(Size of the dispensable genome)]\label{T4}
For $G$ as above,
%
%
\begin{equation}
\label{eq:T3a}
\mathbb E[G] = \theta\sum_{i=0}^{n-1} \frac{1}{\rho+i}.
\end{equation}
In addition, with $g_{(k_1,k_2,k_3)}$ given in Definition
\ref{def:gh},
%
%
\begin{equation}
\label{eq:T3b}
\mathbb V[G] = \theta\sum_{i=0}^{n-1} \frac{1}{\rho+i} -
\theta^2 \Biggl(\sum_{i=0}^{n-1} \frac{1}{\rho+i} \Biggr)^2
+ \frac{\theta^2}{4} g_{(n,0,0)}.
\end{equation}
\end{theorem}
\begin{remark}
An estimate for the size of the pangenome (dispensable plus core
genome) in real bacterial populations has attained much interest
[e.g., \citet{pmid16172379}, \citet{Lapierre2009}]. Most
interestingly,
some species like \textit{Bacillus anthracis} seem to have a
\textit{closed genome}, that is, only a limited number of genes in the
pangenome, since no new gene was found after sequencing the fourth
out of eight strains [\citet{pmid16172379}]. Other species like
\textit{Prochlorococcus} and \textit{Synechococcus}, have an \textit{open
genome} since estimates based on 22 different strains show that
every newly sequenced genome exhibits 277 new genes on average
(Baumdicker, unpublished observation).
For open genomes, a model
based on some linguistic insights, Heap's law, has been considered
[\citet{pmid19086349}]. As a result, a power law for the total number
of genes is found and it is estimated that a total of $n^{0.43\pm
0.02}$ genes are found in a sample of $n$ individuals from
\textit{Bacillus cereus}. This finding is in stark contrast to our
theorem, which implies that the number of genes grows only
logarithmically in $n$. However, in order to decide which is the
correct asymptotics certainly requires a lot more data, since $n=14$
strains are not enough to specify asymptotic behavior.

We conjecture that $\mathbb V[G]$ grows like $\mathbb
E[G]$ for $n\to\infty$. (The corresponding statement is true in the
infinite sites model [see \citet{Wakeley2008}, (4.8)].) The
reason is
that for given $\mathcal T$, $G$ is Poisson distributed with a
parameter increasing with the tree length. In addition, for large
$n$ the length of the Kingman coalescent is largest near the leaves
and the coalescent almost becomes deterministic near the leaves.
For example, it has been shown that the sum of \textit{external branch}
lengths (i.e., branches connecting a leaf to the next node in the
tree) converges to 2 in $L^2$ [\citet{Fu1995}].

\setcounter{footnote}{1}

We give an example for the computation of
$g_{\underline k}$ in the case $\underline k = (2,0,0)$. For the
calculation, we observe that $\lambda_1=1,
\lambda_2=\lambda_3=\lambda_4 = 0, \lambda_5=\lambda_6=\rho$,
$\overline\lambda= 1+2\rho$ and, from (\ref{eq:defgh0}) and
(\ref{eq:defgh1}),
\begin{eqnarray*}
h_1 & = & \frac2\rho,\qquad h_2 = \frac{2}{\rho} +
\frac{2}{1+\rho} = \frac{2(1+2\rho)}{\rho(1+\rho)},\\
g_{(1,0,0)} & = & \frac{4}{\rho^2},\qquad g_{(1,1,0)} =
g_{(1,0,1)} = \frac4 \rho\biggl(\frac1 \rho+
\frac{1}{\rho+1} \biggr) = \frac{4(1+2\rho)}{\rho^2(1+\rho)}.
\end{eqnarray*}
The recursion (\ref{eq:defgh2}) then gives
%
%
\begin{eqnarray}
\label{eq:g200}\qquad
g_{(2,0,0)} &=& \frac{8}{(1+2\rho)^2} + \frac{1}{1+2\rho} \biggl(
\frac{4h_1}{1+2\rho} + g_{(1,0,0)} \biggr)  \nonumber\\
&&{}+\frac{2\rho}{1+2\rho} \biggl(
\frac{1}{1+2\rho}(2h_1+2h_2) + g_{(1,0,1)} \biggr) \nonumber\\
&=& 4 \biggl(\frac{2}{(1+2\rho)^2} + \frac{2}{\rho(1+2\rho)^2} +
\frac{1}{\rho^2(1+2\rho)} + \frac{2}{(1+2\rho)^2} \nonumber\\
&&\hspace*{105.5pt}{}
+
\frac{2}{(1+2\rho)(1+\rho)} + \frac{2}{\rho(1+\rho)} \biggr)
\\
&=&
4 \biggl( \frac{1}{\rho^2} + \frac{2}{\rho(1+\rho)} +
\frac{2}{(1+2\rho)(1+\rho)} \biggr) \nonumber\\
&=& 4 \biggl(\frac{1}{\rho} +
\frac{1}{(1+\rho)} \biggr)^2 + \frac{4}{(1+\rho)^2(1+2\rho)}.\nonumber
\end{eqnarray}
Using (\ref{eq:T3b}) this then gives for $n=2$
\[
\mathbb V_{n=2}[G] = \theta\frac{1+2\rho}{\rho(1+\rho)} +
\theta^2 \frac{1}{(1+\rho)^2(1+2\rho)}.
\]
For $n=3$, the computation is more involved\footnote{Several
computations in the paper are most easily done using a program
like \textsc{Mathematica}. Therefore, a \textsc{Mathematica}-notebook
with all relevant computations can be downloaded from the homepage
of the corresponding author.} and leads to
\[
\mathbb V_{n=3}[G] = \frac{\theta}{\rho} + \frac{\theta}{1+\rho}
+ \frac{\theta}{2+\rho} + \theta^2 \frac{90 + 249 \rho+ 275
\rho^2 + 145 \rho^3 + 30 \rho^4}{(1 + \rho)^2 (2 + \rho)^2 (1 +
2 \rho) (3 + 2 \rho) (6 + 5 \rho)}.
\]
\end{remark}

\subsection{\texorpdfstring{Gene frequency
spectrum.\setcounter{footnote}{2}\protect\footnote{The term \textit{gene
frequency spectrum} was
used by \citet{Kimura1964} to denote the frequency of alleles in
the infinite sites model. Later, the term changed to \textit{site
frequency spectrum} since single sites on the chromosome could
be sequenced [e.g., \citet{Durrett2008}]. Here, we reintroduce
the term for gene frequencies in the infinitely many genes
model.}}{Gene frequency spectrum}}
\label{S:35}
By definition, core genes are present in all individuals of the
sample. In contrast, genes from the dispensable genome can be present
at any frequency. These possibilities give rise to the gene frequency
spectrum.

The \textit{gene frequency spectrum} (\textit{of the dispensable
genome}) is
given by $G_1,\ldots, G_n$, where
%
%
\begin{equation}\label{eq:Gk}
G_k^{(n)} := G_k := |\{u\in I\dvtx u\in\mathcal G_i \mbox{ for exactly }
k \mbox{ different }i\}|.
\end{equation}

\begin{theorem}[(Gene frequency spectrum)]\label{T5}
For $G_1,\ldots,G_n$ as above,
\[
\mathbb E[G_k] = \frac{\theta}{k}\frac{n\cdots
(n-k+1)}{(n-1+\rho)\cdots(n-k+\rho)},\qquad k=1,\ldots,n.
\]
\end{theorem}
\begin{remark}
In the case $\rho=0$, genes cannot get lost and consequently
$G_n=\infty$. However, the classes $k=1,\ldots,n-1$ consist of genes
segregating in the sample (since both states---presence and
absence of the gene---are observed). Hence, as discussed in
Section \ref{S:23}, these classes follow predictions for the
infinite sites model. In this model, it is implicit in results
already obtained by \citet{Wright1938} [and later were refined by
\citet{Kimura1964}, \citet{Griffiths2003}, \citet
{EvansShvetsSlatkin2007}] that
\[
\mathbb E_{\rho=0}[G_k] = \frac\theta k.
\]
On the other hand, by the theorem,
\[
\mathbb E_\rho[G_k]
\stackrel{\rho\downarrow0}{\hbox to 1cm{\rightarrowfill}}
\frac{\theta n}{k(n-k)}
\]
such that
the gene frequency spectrum is not continuous at $\rho=0$.

The model for the bacterial supragenome,
introduced in \citet{pmid16172379} takes population frequencies of
genes into account, that is, the gene frequency spectrum. While the
supragenome model assumes several different frequency classes to
begin with, we derive the gene frequency spectrum from first
principles, that is, from gene gain and loss events along the
genealogy.
\end{remark}

\subsection{\texorpdfstring{Union of core and dispensable
genome.}{Union of core and dispensable genome}}
\label{S:36}
Until now we only derived results for the dispensable genome. In data
obtained from bacterial species, the union of the core and
dispensable
genome is of primary interest. It is straightforward to extend our
results to this union:
\eject

If we replace $\mathcal G_i$ by $\mathcal G_i\cup\mathcal G_c$, $1\leq
i\leq n$, in (\ref{eq:Hbar})--(\ref{eq:Aijkl}),
(\ref{eq:G}) and (\ref{eq:Gk}), recall $g_c:=|\mathcal G_c|$, and
denote the resulting quantities by $\widetilde A$, $\widetilde D$,
$\widetilde P$, $\widetilde D_{ij,kl}$, $\widetilde G$, $\widetilde
G_k$, we obtain
\[
\widetilde A = A + g_c,\qquad \widetilde D = D,\qquad \widetilde
D_{ij,kl} = D_{ij,kl},\qquad \widetilde P = P,\qquad \widetilde G
= G + g_c
\]
and
\[
\widetilde G_k = \cases{ G_k, &\quad $k=1,\ldots,n-1$, \cr
G_k + g_c, &\quad $k=n$.}
\]
Hence, properties of $\widetilde A, \widetilde D, \widetilde P,
\widetilde G, \widetilde G_k$ follow immediately from Theorems
\ref{T1}--\ref{T5}.

\subsection{\texorpdfstring{Application: A dataset from
\textup{Prochlorococcus}.}{Application: A dataset from
\textup{Prochlorococcus}}}
\label{S:37}
Data from complete gen\-omes of a population sample of bacteria have
been available only for a few years. Because the infinitely many genes
model we propose is new in the population genetic context, we show
some data in order to see if the model as studied above could be
realistic.

Here we chose a set of $n=9$ strains of \textit{Prochlorococcus} which
appear to be closely related. \textit{Prochlorococcus} is a marine
picocyanobacterium (length $\sim$ $0.6$~$\mu$m, genome size $\sim2$ Mbp)
living\vspace*{1pt} in the ocean at depth up to 200 m. Their population size can be
as large as $10^6$ individuals (i.e., cells) per ml. In total, 22
complete genome sequences of these cyanobacteria are available in
GenBank at the moment [\citet{pmid18159947},
\citet{Dufresne2008GenomeBiol18507822}]. The $n=9$ chosen
\textit{Prochlorococcus} genomes are similar to each other in terms of
GC-content and share a similar physiology.

\begin{figure}

\includegraphics{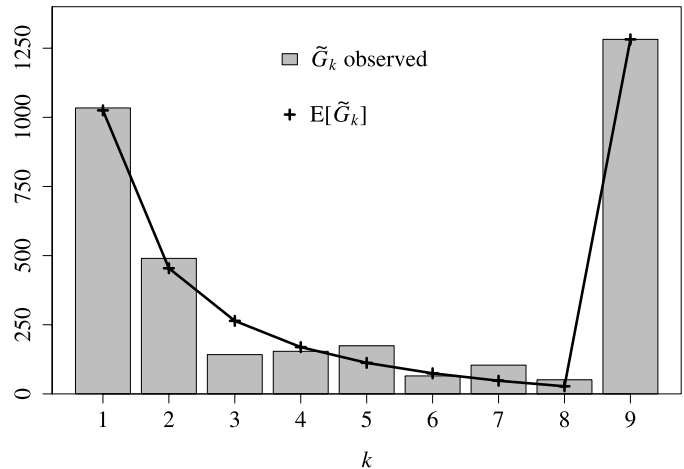}

\caption{The fit of observed data from nine
closely related strains of \textup{Prochlorococcus} with the
expectations for the gene frequency spectrum. Estimates were as
given in (\protect\ref{eq:est}).}\label{fig:freqSpec}
\end{figure}

We estimate the model parameters $\theta, \rho$ and $g_c$ based on the
gene frequency spectrum $\widetilde G_1,\ldots, \widetilde G_9$ which we
compare with our results from Theorem \ref{T5}. The number of genes
present in all individuals is 1282, forming the largest class in the
observed gene frequency spectrum (see Figure
\ref{fig:freqSpec}). Genes occurring in only a single individual were
the second largest class with 1034 genes. By a least squares fit of
$\widetilde G_k$ and $\mathbb E[\widetilde G_k]$ for $k=1,\ldots,n$ we
obtain the estimates
%
%
\begin{equation}\label{eq:est}
\widehat\theta= 1142.17,\qquad \widehat\rho= 2.03,\qquad
\widehat g_c = 1270.
\end{equation}
Note that the estimate for $g_c$ means that we expect that 14 genes
which are carried by all individuals belong to the dispensable
genome. As shown in Figure \ref{fig:freqSpec}, these estimates produce
%
%
a reasonably good fit with the data. Of course, a statistical test
which is able to reject our model for gene content in general, and the
assumption that all genes in the dispensable genome evolve neutrally
in particular, would be desirable.

\subsection{\texorpdfstring{Outlook.}{Outlook}}
\label{S:out}
We introduce the infinitely many genes model on a Kingman coalescent
as a simple null-model of genome evolution in bacterial
species. However, both the reproduction and the mutation dynamics can
be extended to become biologically more realistic. For the
reproduction dynamics, several extensions have been considered in the
literature, for example, structured populations and populations of varying
size [see, e.g., \citet{Durrett2008}].

The mutation dynamics can be extended as well. Our strongest
assumption is that genes taken from the environment are completely
new. In particular, the model does not allow for genes being
transferred between individuals directly. Such a physical exchange of
genes between bacteria is known as horizontal gene transfer. The
underlying mechanism is \textit{bacterial conjugation}. The donor cell
produces a pilus that attaches to the recipient cell and a
single strand of DNA is transported from the donor to the
recipient. After replication of the DNA, both cells carry the
transferred genetic material. The duration of conjugation is long
enough in order to transfer several genes. Hence, by events of
horizontal gene transfer, the transferred genes do not share the
genealogy of the cell line. Thus, building such a mechanism into the
above model requires the use of different genealogical trees for
different genes. Such a mechanism was already considered in the
phylogenetics literature by \citet{pmid12874054}.

In order to add even more biological realism, at least three aspects
can be considered:
\begin{enumerate}
\item As \citet{pmid17475002} show there is frequent
recombination even
within the core genome. Such recombination can also be explained by
conjugation and has attained much interest [e.g.,
\citet{pmid17255503}] since the amount of recombination is known
to be related to sequence similarity [e.g., \citet{pmid9275198}],
suggesting that bacterial species can be distinguished by the extent
of recombination between strains [\citet{DykhuizenGreen1991},
\citet{MaynardSmith1995}].
\item As seen in genomic data, several genes are clustered in gene
families. This is best explained by events of \textit{gene
duplication} with a potential \textit{subfunctionalization} of these
genes along ancestral lines [e.g., \citet{Durrett2008},
\citet{DurrettPopovic2007}].
\item There are certainly selective constraints on the number of genes
in the dispensable genome. If these genes are evolving neutrally and
are not necessarily needed for a bacterium to function properly,
selection should act in order to minimize the dispensable genome.
\end{enumerate}

Specifying the set of genes of an individual requires that the whole
genome of the individual is sequenced. Finding the different genes in
a dataset like the one used in the last section means that open
reading frames (ORFs), that is, regions in the genome between start
codons and stop codons of all individuals are found. In the dataset we
say that two individuals carry the same gene if we find a pair of ORFs
in both individuals that are highly similar. However, the DNA
sequence of this pair of ORFs is usually not identical. Refined
mutation models should extend our approach and describe the genomic
diversity of the different genes as well as the variation of DNA
sequences within the genes.

\section{\texorpdfstring{The one-line-equilibrium and proof of Theorem
\protect\ref{T1}.}{The one-line-equilibrium and proof of Theorem 1}}
\label{S:proof1}
Consider a sample of size $n$ and recall the sets of genes $\mathcal
G_1,\ldots,\mathcal G_n$ from Section \ref{S:22}. All results we provide
with Theorems \ref{T1}--\ref{T5} are dealing with the joint
distribution of $\mathcal G_1,\ldots,\mathcal G_n$. We start with
properties of one- and two-dimensional marginals of the total masses
of this joint distribution. First, we have to obtain a key result for
the gene content along a single ancestral line in Section
\ref{S:41}. The first two moments of the one- and two-dimensional
marginals are obtained in Section \ref{S:42} which then lead to a
proof of Theorem \ref{T1} in Section \ref{S:43}.

\subsection{\texorpdfstring{The one-line equilibrium.}{The one-line
equilibrium}}
\label{S:41}
We start with some arguments that will appear frequently in the next
sections. For $n=1$, the random tree $\mathcal T$ is only a single
infinite line. We consider the gene content along a single ancestral
line $\mathcal T = \mathbb R_-$. In this setting, recall the process
$\Gamma_{\mathbb R_-} = (\mathcal G_t)_{t\in\mathbb R_-}$ from
Definition \ref{def:21}. Note that, almost surely, $\mathcal G_t$ does
not have double points for all $t\in\mathbb R_-$. Recall our notation
from Remark \ref{rem:finMeas}.
\begin{definition}[(Poisson random measure and thinning)]
We denote by
$\mathcal{POI}(\alpha)$ the distribution of a
Poisson random measure with intensity measure~$\alpha$. We will
also write $\mathcal{POI}(\alpha)$ for the Poisson distribution
with parameter $\alpha$ if $\alpha\in\mathbb R_+$.

For $g\in\mathcal N_f(I)$, we denote by $\mathcal{THIN}(g,p)$
the distribution of the random measure arising by keeping any
point in $g$ with probability $p$.
\end{definition}
\begin{proposition}[(Distribution of $\Gamma_{\mathbb R_-}$)]
\label{P:41}
Let $s<t$.

Given $\mathcal G_s=g \in\mathcal N_f(I)$, the two random
measures $\mathcal G_t \cap\mathcal G_s$ and $\mathcal G_t
\setminus\mathcal G_s$ are independent. Their distribution is
given by
\begin{eqnarray*}
\mathcal G_t\cap\mathcal G_s
& \sim& \mathcal{THIN}\bigl(g,e^{-\rho/2 (t-s)}\bigr), \\
\mathcal G_t \setminus\mathcal G_s & \sim&
\mathcal{POI}\biggl(\frac\theta\rho\bigl(1-e^{-\rho/2
(t-s)}\bigr)\cdot\lambda_I\biggr).
\end{eqnarray*}
%
%

The distribution
$\mathcal{POI}(\frac\theta\rho\lambda_I)$ is the unique
equilibrium for $\Gamma_{\mathbb R_-}$ and it is reversible. In
equilibrium, $\mathcal G_t\cap\mathcal G_s, \mathcal G_t
\setminus\mathcal G_s$ and $\mathcal G_s \setminus\mathcal G_t$
are independent and their distributions are given by
\begin{eqnarray*}
\mathcal G_t\cap\mathcal G_s & \sim&
\mathcal{POI}\biggl(\frac\theta\rho e^{-\rho/2
(t-s)}\cdot\lambda_I\biggr), \\
\mathcal G_s \setminus\mathcal G_t
&\stackrel{d}{=}& \mathcal G_t \setminus\mathcal G_s \sim
\mathcal{POI}\biggl(\frac\theta\rho\bigl(1-e^{-\rho/2
(t-s)}\bigr)\cdot\lambda_I\biggr).
\end{eqnarray*}
\end{proposition}
\begin{remark}
Recall from Remark \ref{rem:finMeas} that we identify $\mathcal G_s,
\mathcal G_t$ with the set of genes carried at times $s$ and
$t$. Note that $\mathcal G_s\cap\mathcal G_t = \mathcal G_s\wedge
\mathcal G_t$ represents the genes~present both at time $s$ and at
time $t$. Moreover, $\mathcal G_s\setminus\mathcal G_t = (\mathcal
G_s-\mathcal G_t)^+$ are the genes present at time $s$ but absent at
time $t$, that is, genes lost during time $(s,t]$. The genes in
$\mathcal G_t\setminus\mathcal G_s = (\mathcal G_t-\mathcal G_s)^+$
are genes gained during time $(s,t]$. As the proposition shows, all
three quantities are independent in equilibrium.
\end{remark}
\begin{pf*}{Proof of Proposition \ref{P:41}}
First, recall that all new points in $\mathcal G_{s'},
s<s'\leq t$ are pairwise different and different from points in
$\mathcal G_s$, almost surely. During $(s,t]$ several points of
$\mathcal G_{s'}$ are lost. A point in $\mathcal G_{s'}$ is not lost
with probability $e^{-\rho/2 (t-s')}$. Since all points are
lost independently, we find that $\mathcal G_s \cap\mathcal G_t
\sim\mathcal{THIN}(g,e^{-\rho/2 (t-s)})$. Additionally,
several new points in $\mathcal G$ arise during $(s,t]$. Hence, we
find that $\mathcal G_t \setminus\mathcal G_s$ is independent of
both, $\mathcal G_s$ and $\mathcal G_s\cap\mathcal G_t$. To obtain
the distribution of $\mathcal G_t \setminus\mathcal G_s$, note that
a point in $\mathcal G_{s'}\setminus\mathcal G_{s'-}$ is lost at
rate $\frac\rho2 $ and hence is present in $\mathcal G_t$ with
probability $e^{-\rho/2 (t-s')}$. Since new\vspace*{1pt} points arise at
rate $\frac\theta2$ during $(s,t]$ and are lost independently, we
find that\vspace*{1pt} the number of points in $\mathcal G_t \setminus\mathcal
G_s$ is Poisson distributed with parameter $\frac\theta
2\int_0^{t-s} e^{-\rho/2 (t-s')}\,ds' =
\frac\theta\rho(1-e^{-\rho/2 (t-s)})$. Since these points
must be uniformly distributed on $I$, we have that $\mathcal G_t
\setminus\mathcal G_s \sim
\mathcal{POI}(\frac\theta\rho(1-e^{-\rho/2
(t-s)})\cdot\lambda_I)$. So we have shown the first assertion.

To see that $\mathcal{POI}(\frac\theta\rho\cdot
\lambda_I)$ is the unique equilibrium of $\Gamma_{\mathbb R}$, note
that there can be at most one equilibrium since the Markov process
$\Gamma_{\mathbb R_-}$ is Harris recurrent. Moreover, if $\mathcal
G_s \sim\mathcal{POI}(\frac\theta\rho\cdot\lambda_I)$, then
$\mathcal{THIN}(\mathcal G_s, e^{-\rho/2 (t-s)}) =
\mathcal{POI}(\frac\theta\rho\times\break e^{-\rho/2
(t-s)}\cdot\lambda_I)$ and so
\[
\mathcal G_t \sim
\mathcal{POI}\biggl(\frac\theta\rho e^{-\rho/2
(t-s)}\cdot\lambda_I\biggr) \ast
\mathcal{POI}\biggl(\frac\theta\rho\bigl(1-e^{-\rho/2
(t-s)}\bigr)\cdot\lambda_I\biggr) =
\mathcal{POI}\biggl(\frac\theta\rho\cdot\lambda_I\biggr)
\]
(where $\ast$
denotes convolution). For reversibility, we write $\mathcal G_s =
(\mathcal G_s\cap\mathcal G_t) \uplus(\mathcal G_s
\setminus\mathcal G_t), \mathcal G_t = (\mathcal G_s\cap\mathcal
G_t) \uplus(\mathcal G_t \setminus\mathcal G_s)$ where $\mathcal
G_s\cap\mathcal G_t, \mathcal G_s \setminus\mathcal G_t,
\mathcal G_t \setminus\mathcal G_s$ are independent, such that
$\Gamma_{\mathbb R_-}$ is in equilibrium at times $s$ and $t$ and
\begin{eqnarray*}
\mathcal G_s\cap\mathcal G_t & \sim&
\mathcal{POI}\biggl(\frac\theta\rho e^{-\rho/2
(t-s)}\cdot\lambda_I\biggr),\\
\mathcal G_s \setminus\mathcal G_t
&\stackrel{d}{=}& \mathcal G_t \setminus\mathcal G_s
\sim\mathcal{POI}\biggl(\frac\theta\rho\bigl(1-e^{-\rho/2
(t-s)}\bigr)\cdot\lambda_I\biggr).
\end{eqnarray*}
By this representation, given some continuous functions $f_1,f_2\dvtx
I\to\mathbb R$, writing $\langle f_i, x\rangle:= \int f_i \,dx,
i=1,2$,
\begin{eqnarray*}
\mathbb E\bigl[e^{-\langle f_1, \mathcal G_s\rangle} \cdot e^{-\langle
f_2, \mathcal G_t\rangle}\bigr]
& = & \mathbb E\bigl[e^{-\langle f_1+f_2,
\mathcal G_s\cap\mathcal G_t\rangle}\bigr] \cdot\mathbb
E\bigl[e^{-\langle f_1, \mathcal G_s\setminus\mathcal G_t\rangle}\bigr]
\cdot\mathbb E\bigl[e^{-\langle f_2, \mathcal G_t\setminus\mathcal
G_s\rangle}\bigr] \\
& = & \mathbb E\bigl[e^{-\langle f_1+f_2, \mathcal
G_s\cap\mathcal G_t\rangle}\bigr] \cdot\mathbb E\bigl[e^{-\langle f_1,
\mathcal G_t\setminus\mathcal G_s\rangle}\bigr] \cdot\mathbb
E\bigl[e^{-\langle f_2, \mathcal G_s\setminus\mathcal G_t\rangle}\bigr] \\
& = & \mathbb E\bigl[e^{-\langle f_2, \mathcal G_s\rangle} \cdot
e^{-\langle f_1, \mathcal G_t\rangle}\bigr].
\end{eqnarray*}
Hence, since the joint Laplace transforms $ \mathbb E[e^{-\langle
f_1, \mathcal G_s\rangle} \cdot e^{-\langle f_2, \mathcal
G_t\rangle}]$ determine the joint distribution of $(\mathcal G_s,
\mathcal G_t)$ uniquely, we find that $(\mathcal G_s, \mathcal G_t)
\stackrel{d}{=} (\mathcal G_t, \mathcal G_s)$ and reversibility is
shown.
\end{pf*}

\subsection{\texorpdfstring{Gene content in individuals and
pairs.}{Gene content in individuals and pairs}}
\label{S:42}
Next we obtain the first two moments of the two-dimensional
distribution of $(|\mathcal G_1|, \ldots, |\mathcal G_n|)$.
\begin{proposition}[{[Distribution of $(\mathcal G_i, \mathcal G_j)$]}]
\label{P:42}
For $i=1,\ldots,n$,
\[
\mathcal G_i \sim\mathcal{POI} \biggl(\frac\theta\rho\cdot\lambda_I
\biggr).\vadjust{\goodbreak}
\]

In particular,
\[
\mathbb E[|\mathcal G_i|] = \mathbb V[|\mathcal G_i|] = \frac\theta
\rho.
\]

For $1\leq i\neq j\leq n$,
\[
\mathbb{COV}[|\mathcal G_i|, |\mathcal G_j|] = \frac{\theta}{\rho
(1+\rho)}.
\]
\end{proposition}
\begin{remark}
In the proof of 2. we use the well-known fact that for random
variables $X, Y, T$
\[
\mathbb{COV}[X,Y] = \mathbb{COV} [\mathbb E[X|T], \mathbb E[Y|T] ]
+ \mathbb E [\mathbb{COV}[X,Y|T] ]
\]
with
\[
\mathbb{COV}[X,Y|T] := \mathbb E [(X-\mathbb E[X|T])
(Y-\mathbb E[Y|T]) |T ].
\]
\end{remark}
\begin{pf*}{Proof of Proposition \ref{P:42}}
Consider the ancestral line of individual $i$. The
process $(\mathcal G_t)_{t\preceq i}$ follows the same dynamics as
the process $\Gamma_{\mathbb R_-}$ studied in Section~\ref{S:41}.
The first claim follows from Proposition \ref{P:41}(2), which shows
that the unique reversible equilibrium for this dynamics is
$\mathcal{POI}(\frac\theta\rho\lambda_I)$.

We denote the random coalescence time of individuals
$i$ and $j$ by $T$. As above, $\mathcal G_i$ and $\mathcal G_j$
denote (the finite measures describing) the genes present in both
individuals. Recall that we have shown in Proposition \ref{P:41}(2)
that the equilibrium of the Markov chain $\Gamma_{\mathbb R_-} =
(\mathcal G_t)_{t\in\mathbb R_-}$ of Section \ref{S:41} is
reversible. Hence, given $T$, we have that $(\mathcal G_i, \mathcal
G_j)$ and $(\mathcal G_{-2T}, \mathcal G_{0})$ have the same
distribution. So we find that $\mathcal G_i \cap\mathcal G_j,
\mathcal G_i \setminus\mathcal G_j$ and $\mathcal G_j \setminus
\mathcal G_i$ are independent and, by Proposition \ref{P:41},
\[
\mathcal G_i \cap\mathcal G_j \sim\mathcal{POI}\biggl(\frac\theta\rho
e^{-\rho T}\cdot\lambda_I\biggr),\qquad
\mathcal G_i \setminus\mathcal G_j \stackrel{d}{=} \mathcal
G_i \setminus\mathcal G_j \sim
\mathcal{POI}\biggl(\frac\theta\rho(1-e^{-\rho T})\cdot\lambda_I\biggr).
\]
Moreover, both sets of genes, $\mathcal G_i$ and $\mathcal G_j$, are
independent of $T$. We obtain
\begin{eqnarray*}
&&
\mathbb{COV} [|\mathcal G_i|,
|\mathcal G_j| ] \\
&&\qquad = \mathbb{COV} [ \mathbb E[|
\mathcal G_i| |T ],\mathbb E[|\mathcal G_j| | T] ]
+ \mathbb E [ \mathbb{COV}[|\mathcal G_i|,
|\mathcal G_j| | T] ]\\
&&\qquad = \mathbb{COV} \biggl[ \frac{\theta}{\rho},\frac{\theta}{\rho} \biggr]
+ \mathbb E \bigl[ \mathbb{COV}[|\mathcal G_i\cap\mathcal G_j| +
|\mathcal G_i \setminus\mathcal G_j| , |\mathcal G_i\cap\mathcal
G_j| +
|\mathcal G_j \setminus\mathcal G_i| | T] \bigr]\\
&&\qquad = \mathbb E [ \mathbb V[|\mathcal G_i\cap\mathcal G_j| |T] ]
= \mathbb E \biggl[ \frac{\theta}{\rho}
e^{- \rho T} \biggr] = \frac{\theta}{\rho(1+\rho)}
\end{eqnarray*}
as $T\sim\operatorname{Exp}(1)$.
\end{pf*}

\subsection{\texorpdfstring{Proof of Theorem \protect\ref{T1}.}{Proof
of Theorem 1}}
\label{S:43}
Theorem \ref{T1} now follows from Proposition \ref{P:42} and
\begin{eqnarray*}
\mathbb E[A] & = & \frac{1}{n}\sum_{i=1}^n \mathbb E[|\mathcal G_i|] =
\frac\theta\rho,\\
\mathbb V[A] & = & \frac{1}{n^2} \Biggl(\sum_{i=1}^n \mathbb V[|\mathcal
G_i|] + \mathop{\sum_{i,j=1}}_{i\neq j}^n
\mathbb{COV}[|\mathcal G_i|, |\mathcal G_j|] \Biggr) \\
& = & \frac1n
\frac\theta\rho+ \biggl( 1 - \frac1n \biggr)
\frac{\theta}{\rho(1+\rho)} = \frac1n \frac{\theta}{1+\rho} +
\frac{\theta}{\rho(1+\rho)}.
\end{eqnarray*}
%

\section{\texorpdfstring{Extension of Proposition
\protect\ref{P:41} and proof of Theorem \protect\ref{T2}.}{Extension of
Proposition 4.2 and proof of Theorem 2}}
The one-line equilibrium considered in Proposition \ref{P:41} provides
the right setting for computing the one- and two-dimensional marginals
of $\mathcal G_1, \ldots, \mathcal G_n$ as shown in the proof of
Proposition \ref{P:42}. In Section \ref{S:51} we provide a method to
compute higher order marginals. We will use this method for second
(Section \ref{S:52}), third (Section \ref{S:53}) and fourth (Section
\ref{S:54}) order which finally leads to a proof of Theorem \ref{T2}
in Section \ref{S:55}.

\subsection{\texorpdfstring{Extending the one-line equilibrium to a genealogical
tree.}{Extending the one-line equilibrium to a genealogical tree}}
\label{S:51}
Before we introduce the general method, how to obtain all marginals of
$\mathcal G_1, \ldots, \mathcal G_n$, we have to set the scene. Consider
the genealogical tree $\mathcal T$ relating all $n$ individuals and
the tree-indexed Markov chain $\Gamma_{\mathcal T} = (\mathcal
G_t)_{t\in\mathcal T}$. In equilibrium, we have seen above that
$\mathcal G_t \sim\mathcal{POI}(\frac\theta\rho\cdot\lambda_I)$ for
all $t\in\mathcal T$. Hence, we now consider the case that $\mathcal
T$ is a rooted tree with root $r$ and $\mathcal G_r \sim
\mathcal{POI}(\frac\theta\rho\cdot\lambda_I)$. We need some notation
to deal with the genealogical tree $\mathcal T$.
%
%
\begin{definition}[(Survival function)]
\label{def:tree}
Let $\mathcal T$ be a binary tree with one distinguished point
$r\in\mathcal T$, referred to as the root of $\mathcal T$, a finite
set of leaves $\mathcal L \subseteq\mathcal T$ and internal
vertices $\mathcal V$. For $s,t\in\mathcal T$ we denote by $(s,t]$
the set of points which must be visited on any path between $s$ and
$t$. Moreover, $d_{\mathcal T}(s,t)$ is the length of the path
between $s$ and $t$.
Define a partial order $\preceq$ on $\mathcal T$ by saying that
$s\preceq t$ iff $s \in(r,t]$ (such that $r$ is the minimal
element). For $s,t\in\mathcal T$ the point $s\wedge t$ is given as
the maximal element in $\{q\dvtx q\preceq s \mbox{ and }q\preceq
t\}$. For an internal node (i.e., a branch point) $t\in\mathcal T$
we denote by $t_{1}$ and $t_{2}$ the two directions in $\mathcal T$
leading to bigger (with respect to $\preceq$) elements.

We define the \textit{survival function} $p_{\mathcal T}\dvtx\mathcal
T \to[0,1]$ by
%
%
\begin{eqnarray}
\label{eq:52}
p_{\mathcal T}(t) & = & 1 \qquad\mbox{for } t\in\mathcal L,\nonumber\\
\frac{\partial p_{\mathcal T}(t)}{\partial t} & = & \frac\rho2
p_{\mathcal T}(t) \qquad\mbox{for } t\in\mathcal T \setminus(\mathcal L
\cup\mathcal V),\\
p_{\mathcal T}(t) & = & 1 - \bigl(1-p_{\mathcal T}(t_{1})\bigr)\bigl(1-p_{\mathcal
T}(t_{2})\bigr) \qquad\mbox{for } t \in\mathcal V,\nonumber
\end{eqnarray}
where for $f\dvtx\mathcal T\setminus(\mathcal L\cup\mathcal V)\to\mathbb
R$
\[
\frac{\partial f(t)}{\partial t} := \lim_{\varepsilon\to0} \frac1
\varepsilon\bigl(f(t+\varepsilon) - f(t)\bigr)
\]
and $t+\varepsilon$ is any
point in $\mathcal T$ with $d_{\mathcal
T}(t,t+\varepsilon)=\varepsilon$ and $t\preceq t+\varepsilon$,
if the limit exists.
\end{definition}
\begin{proposition}[(Probability of no loss along $\mathcal
T$)]\label{P:52}
Let $\mathcal T$ be a binary tree, rooted at $r$, $p_{\mathcal T}$
as in Definition \ref{def:tree} and $\Gamma_{\mathcal T} = (\mathcal
G_t)_{t\in\mathcal T}$ be the tree-indexed Markov chain from
Definition \ref{def:21} with $\mathcal G_r \sim
\mathcal{POI}(\frac\theta\rho\lambda_I)$. Then for $u\in I$ and
$t\in\mathcal T$
\[
\mathbb P \Biggl[ u\in\mathop{\bigcup_{t\preceq s}}_{s\in\mathcal L}
\mathcal G_s \Big| u\in\mathcal G_t \Biggr] = p_{\mathcal T}(t).
\]
\end{proposition}
\begin{pf}
Denote the probability on the left-hand side by $q(t)$. First note
that $q(t)=1$ if $t\in\mathcal L$ since $\{s\in\mathcal L\dvtx t\preceq
s\} = \{t\}$. Moreover, the probability on the left-hand side
decreases exponentially at rate $\frac\rho2$ along branches of
$\mathcal T$ due to loss events of $u$. Last, consider the case
$t\in\mathcal V$. Then, $u$ must not be lost to either $t_{1}$ or
$t_{2}$. This occurs with probability $q(t) =
1-(1-q(t_{1}))(1-q(t_{2}))$. In other words, the function $q$
fulfills all defining properties of $p_{\mathcal T}$ from
(\ref{eq:52}) and we are done.
\end{pf}

We need some more notation for subsets of a finite binary rooted tree
$\mathcal T$.
\begin{definition}[(Length and subtrees of $\mathcal
T$)]\label{def:tree2}
We use the notation of Definition \ref{def:tree}.

For the binary tree $\mathcal T$ we denote by $\ell(\mathcal
T)$ its total length, that is, the sum of lengths of all its
branches.

Let $\mathcal L', \mathcal M'\subseteq\mathcal L$ be
sets of leaves with $\mathcal L'\cap\mathcal M'=\varnothing$. We set
$r_0:=\bigwedge_{t\in\mathcal L'} t$ and denote by $\mathcal
T^0(\mathcal L')$ the minimal connected, binary tree spanning the
leaves~$\mathcal L'$, rooted at $r_0$. The set $(\mathcal
T^0(\mathcal L'\cup\mathcal M'))\setminus\mathcal T^0(\mathcal
L')$ consists of $k\leq|\mathcal M'|$ different connected
subtrees, connected with $\mathcal T^0(\mathcal L')$ at vertices
$r_1,\ldots,r_k$. We denote the resulting binary trees by $\mathcal
T^1(\mathcal L', \mathcal M'),\ldots, \mathcal T^k(\mathcal L',
\mathcal M')$, rooted at $r_1,\ldots,r_k$, respectively.
\end{definition}
\begin{remark}\label{rem:onepoint}
For an illustration of the objects introduced in
Definition~\ref{def:tree2}, see Figure \ref{fig:prop}.

%
\begin{figure}

\includegraphics{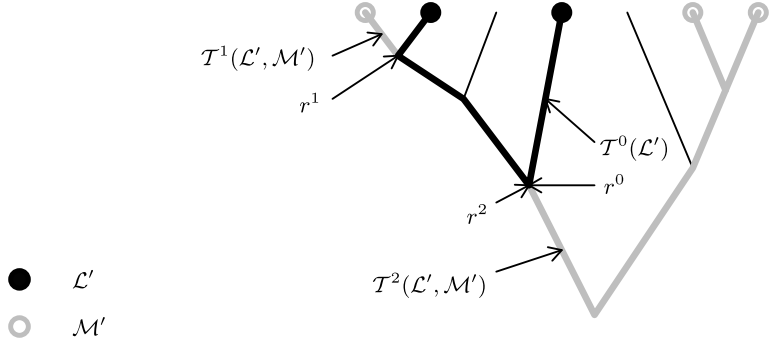}

\caption{Illustration of concepts given in Definition
\protect\ref{def:tree2}(2). The subtree $\mathcal T^0(\mathcal L')$ is
spanned by leaves in $\mathcal L'$. Considering $\mathcal T\setminus
\mathcal T^0(\mathcal L')$, the tree falls in three parts. Two of
them, which lead to a leaf in $\mathcal M'$, are denoted $\mathcal
T^k(\mathcal L', \mathcal M'), k=1,2$. Roots of the trees are $r^0,
r^1, r^2$.}\label{fig:prop}
\end{figure}

If $|\mathcal L'|=1$, it is important to note that $\mathcal
T^0(\mathcal L')$ only consists of a single point. Consequently,
$\ell(\mathcal T^0(\mathcal L'))=0$ in this case.
\end{remark}
\begin{proposition}[(Distribution of $\bigcap\mathcal G_s \setminus
\bigcup\mathcal G_t$)]
\label{P:main}
Let $\mathcal T$ be a finite binary tree, rooted at $r\in\mathcal
T$, $\mathcal L$ its finite set of leaves and $\mathcal L', \mathcal
M'\subseteq\mathcal L$ with $\mathcal L'\cap\mathcal
M'=\varnothing$. Moreover, let $\mathcal T^0(\mathcal L'),\mathcal
T^1(\mathcal L',\mathcal M'),\ldots,\mathcal T^k(\mathcal L',\mathcal
M')$ be as in Definition \ref{def:tree2}. Let $\Gamma_{\mathcal T} =
(\mathcal G_t)_{t\in\mathcal T}$ be the\vspace*{1pt} tree-indexed Markov chain
from Definition \ref{def:21} with $\mathcal G_r \sim
\mathcal{POI}(\frac\theta\rho)$. Then,
\[
\bigcap_{t\in\mathcal L'} \mathcal G_t \Bigm\backslash
\bigcup_{t\in\mathcal M'} \mathcal G_t \sim\mathcal{POI} \Biggl(
\frac\theta\rho e^{-\rho/2 \ell(\mathcal T^0(\mathcal L'))}
\prod_{i=1}^k \bigl(1 - p_{\mathcal T^i(\mathcal L',\mathcal
M')}(r_i)\bigr)\cdot\lambda_I \Biggr).
\]
In addition, if $\mathcal L'', \mathcal M''\subseteq\mathcal L$ with
$\mathcal L'' \cap\mathcal M'' = \varnothing$, then
\[
\bigcap_{t\in\mathcal L'} \mathcal G_t \Bigm\backslash\bigcup_{t\in
\mathcal M'}
\mathcal G_t \quad\mbox{and}\quad \bigcap_{t\in\mathcal L''} \mathcal G_t
\Bigm\backslash\bigcup_{t\in\mathcal M''} \mathcal G_t
\]
are independent if $\mathcal L'\cap\mathcal M''\neq\varnothing$ or
$\mathcal L''\cap\mathcal M'\neq\varnothing$.
\end{proposition}
\begin{remark}
The pairwise independence in the proposition can be extended
to independence of any number of random measures
$\bigcap_{t\in\mathcal L_i} \mathcal G_t \setminus
\bigcup_{t\in\mathcal M_i} \mathcal G_t, i=1,\ldots,n$, provided
$\mathcal L_i\cap\mathcal M_j\neq\varnothing$ or $\mathcal
L_j\cap\mathcal M_i\neq\varnothing$ holds for any pair $i\neq j$.
\end{remark}
\begin{pf*}{Proof of Proposition \ref{P:main}}
Given $\mathcal T$, rooted at $r\in\mathcal T$, we have assumed that
$\mathcal G_r = \mathcal{POI} (\frac\theta\rho
\cdot\lambda_I )$, that is, the tree-indexed Markov chain is in
equilibrium. Consequently, $\mathcal G_{r_0} \sim\mathcal
G_r$. Every gene in $\bigcap_{t\in\mathcal L'} \mathcal G_t
\setminus\bigcup_{t\in\mathcal M'} \mathcal G_t$ must have been
present in $\mathcal G_{r_0}$. In addition, every gene in $\mathcal
G_{r_0}$ has the same chance $p$ to be present in
$\bigcap_{t\in\mathcal L'} \mathcal G_t \setminus
\bigcup_{t\in\mathcal M'} \mathcal G_t$. This\vspace*{1pt} already shows that
$\bigcap_{t\in\mathcal L'} \mathcal G_t \setminus
\bigcup_{t\in\mathcal M'} \mathcal G_t$ is a thinning of a Poisson
measure and hence is Poisson with intensity $\frac\theta\rho p$.
It is important to note that a gene present in
$\bigcap_{t\in\mathcal L'} \mathcal G_t \setminus
\bigcup_{t\in\mathcal M'} \mathcal G_t$ must not be lost on the
whole subtree $\mathcal T^0(\mathcal L')$, which occurs with
probability $e^{-\rho/2 \ell(\mathcal T^0(\mathcal L'))}$ and
must be lost on any subtree leading to a leaf in $\mathcal
M'$. However, the chance that a gene is lost along one such subtree
is given through the survival function. In the subtree $i$, we have
a root $r_i$ connecting the subtree to the
tree spanned by $\mathcal L'$ and so $1 - p_{\mathcal T^i(\mathcal
L',\mathcal M')}(r_i)$ is the probability that the gene is lost in
all leaves in $\mathcal M'$.

For the independence property assume that $\mathcal
L'\cap\mathcal M''\neq\varnothing$ or $\mathcal L''\cap\mathcal
M'\neq\varnothing$. Observe\vspace*{1pt} that $ (\bigcap_{t\in\mathcal L'}
\mathcal G_t \setminus\bigcup_{t\in\mathcal M'} \mathcal G_t )
\cap( \bigcap_{t\in\mathcal L''} \mathcal G_t \setminus
\bigcup_{t\in\mathcal M''} \mathcal G_t ) = 0$ in this case,
that is, $ \bigcap_{t\in\mathcal L'} \mathcal G_t \setminus
\bigcup_{t\in\mathcal M'} \mathcal G_t$ and $\bigcap_{t\in\mathcal
L''} \mathcal G_t \setminus\bigcup_{t\in\mathcal M''} \mathcal
G_t$ arise by different Poisson events along $\mathcal T$. The
independence follows.
\end{pf*}
\begin{corollary}\label{Cor:main}For the same situation as in
Proposition \ref{P:main}, if $\mathcal L'\cap\mathcal M'=\mathcal
L''\cap\mathcal M''=\varnothing$,
\begin{eqnarray*}
&&\mathbb{COV} \biggl[ \biggl|\bigcap_{t\in\mathcal L'}\mathcal G_t
\Bigm\backslash \bigcup_{t\in\mathcal M'} \mathcal G_t \biggr|,
\biggl|\bigcap_{t\in\mathcal L''}\mathcal G_t \Bigm\backslash
\bigcup_{t\in\mathcal M''} \mathcal G_t \biggr| \Big|\mathcal T \biggr]
\\
&&\qquad = \mathbb V \biggl[ \biggl|\bigcap_{t\in\mathcal
L'\cup\mathcal L''}\mathcal G_t \Bigm\backslash\bigcup_{t\in\mathcal
M'\cup\mathcal M''} \mathcal G_t \biggr| \Big|\mathcal T \biggr].
\end{eqnarray*}
\end{corollary}
\begin{pf}
We write
\begin{eqnarray*}
\bigcap_{t\in\mathcal L'}\mathcal G_t \Bigm\backslash
\bigcup_{t\in\mathcal M'} \mathcal G_t
& = &
\biggl(\bigcap_{t\in\mathcal L'}\mathcal G_t \Bigm\backslash
\bigcup_{t\in\mathcal M'\cup\mathcal M''\cup\mathcal L''}\mathcal
G_t \biggr) \uplus\biggl(\bigcap_{t\in\mathcal L'\cup\mathcal
L''}\mathcal G_t \Bigm\backslash\bigcup_{t\in\mathcal M'\cup\mathcal
M''}\mathcal G_t \biggr) \\
&&{} \uplus\biggl(\bigcap_{t\in\mathcal L'\cup
\mathcal L''\cup\mathcal M''}\mathcal G_t \Bigm\backslash
\bigcap_{t\in\mathcal M'}\mathcal G_t
\biggr),\\
\bigcap_{t\in\mathcal L''}\mathcal G_t \Bigm\backslash
\bigcup_{t\in\mathcal M''} \mathcal G_t
& = &
\biggl(\bigcap_{t\in\mathcal L''}\mathcal G_t \Bigm\backslash
\bigcup_{t\in\mathcal M'\cup\mathcal M''\cup\mathcal L'} \mathcal
G_t \biggr) \uplus\biggl(\bigcap_{t\in\mathcal L'\cup\mathcal
L''} \mathcal G_t \Bigm\backslash\bigcup_{t\in\mathcal M'\cup\mathcal
M''}\mathcal G_t \biggr) \\
&&{} \uplus\biggl(\bigcap_{t\in\mathcal L'\cup
\mathcal L''\cup\mathcal M'}\mathcal G_t \Bigm\backslash
\bigcup_{t\in\mathcal M''}\mathcal G_t \biggr).
\end{eqnarray*}
By the independence statement in Proposition \ref{P:main}, only the
covariances of the two second terms in both equalities do not
vanish. The result follows.
\end{pf}

\subsection{\texorpdfstring{Gene content for two individuals.}{Gene
content for two individuals}}
\label{S:52}
The simplest case in Proposition \ref{P:main} arises if $\mathcal T$
has only two leaves. This case was already studied in the proof of
Proposition \ref{P:42}. We extend our analysis by the next result.
\begin{proposition}[(Gene content for two individuals)]
\label{P:56}
For $1\leq i\neq j\leq n$,
\begin{eqnarray*}
\mathbb E[|\mathcal G_i\setminus\mathcal G_j|] & = &
\frac{\theta}{1+\rho},
\\
\mathbb V[|\mathcal G_i\setminus\mathcal G_j|] & = &
\frac{\theta^2}{(1+\rho)^2(1+2\rho)}
+ \frac{\theta}{1+\rho}, \\
\mathbb{COV}[|\mathcal G_i\setminus\mathcal G_j|, |\mathcal
G_j\setminus\mathcal G_i|] & = &
\frac{\theta^2}{(1+\rho)^2(1+2\rho)}.
\end{eqnarray*}
\end{proposition}
\begin{pf}
We use Proposition \ref{P:main}. It suffices to assume that
$\mathcal T$ is a tree connecting individuals $i$ and $j$, that is,
$\mathcal L=\{i,j\}$. First we assume that the coalescence time $T$
of the two individuals is given. Under this assumption, Proposition
\ref{P:main} tells us that
\[
\mathcal G_i \setminus\mathcal G_j \sim
\mathcal{POI}\biggl(\frac\theta\rho(1-e^{-\rho T})\cdot\lambda_I\biggr)
\]
and,
using the fact that expectation and variance are equal for a Poisson
distributed random variable,
\[
\mathbb E [|\mathcal G_i\setminus\mathcal G_j| |T ]
= \mathbb V [|\mathcal G_i\setminus\mathcal G_j| |T ] =
\frac\theta\rho(1-e^{-\rho T}),
\]
such that we obtain
\begin{eqnarray*}
\mathbb E[|\mathcal G_i\setminus\mathcal G_j|] & = & \mathbb E \bigl[
\mathbb E[|\mathcal G_i\setminus\mathcal G_j| | T ] \bigr]=
\mathbb E\biggl[\frac\theta\rho(1-e^{-\rho T})\biggr] = \frac\theta\rho
\frac{\rho}{1+\rho} = \frac{\theta}{1+\rho},\\
\mathbb V[|\mathcal G_i \setminus\mathcal G_j|] & = & \mathbb V
\bigl[ \mathbb E[|\mathcal G_i \setminus\mathcal G_j| |
T] \bigr] + \mathbb E \bigl[ \mathbb V [|\mathcal G_i \setminus
\mathcal G_j| | T] \bigr] \\
& = & \mathbb V
\biggl[\frac{\theta}{\rho}(1-e^{-\rho T}) \biggr] +
\mathbb E \biggl[\frac{\theta}{\rho}(1-e^{-\rho T}) \biggr]\\
& = & \frac{\theta^2}{\rho^2} \biggl( \frac{1}{1+2\rho} -
\frac{1}{(1+\rho)^2} \biggr) +
\frac{\theta}{\rho} \biggl(1-\frac{1}{1+\rho} \biggr) \\
& = &
\frac{\theta^2}{(1+\rho)^2(1+2\rho)} + \frac{\theta}{1+\rho}.
\end{eqnarray*}
In addition, given $T$, $\mathcal G_i\setminus\mathcal G_j$ and
$\mathcal G_j\setminus\mathcal G_i$ are independent by Corollary
\ref{Cor:main}. Hence,
\begin{eqnarray*}
\mathbb{COV}[|\mathcal G_i \setminus\mathcal G_j|, |\mathcal
G_j\setminus\mathcal G_i|]
& = & \mathbb{COV} \bigl[\mathbb E[|\mathcal
G_i\setminus\mathcal G_j| | T],\mathbb E[|\mathcal
G_j\setminus\mathcal G_i| | T] \bigr] \\
& = &\mathbb
V \biggl[\frac\theta\rho(1-e^{-\rho T}) \biggr] =
\frac{\theta^2}{(1+\rho)^2(1+2\rho)}.
\end{eqnarray*}
\upqed\end{pf}

\subsection{\texorpdfstring{Gene content for three individuals.}{Gene
content for three individuals}}
\label{S:53}
Similar to Proposition \ref{P:56}, we use the general setting of
Proposition \ref{P:main} in order to prove results about the joint
distribution of gene content in three individuals.
\begin{proposition}[(Gene content for three individuals)]
\label{P:57}
For $i,j,k\in\{1,\ldots,\break n\}$ pairwise different, 
%
%
\begin{eqnarray}
\label{eq:571}
\mathbb{COV}[|\mathcal G_i \setminus\mathcal G_j|,|\mathcal G_i
\setminus\mathcal G_k|] & = & \frac{\theta^2}{(1+\rho)^2
(1+2\rho)(3+2\rho)} + \frac{\theta}{2+\rho},
\\
\label{eq:572}
\mathbb{COV}[|\mathcal G_i \setminus\mathcal G_j|,|\mathcal G_k
\setminus\mathcal G_i|] & = &
\frac{\theta^2}{(1+\rho)^2(1+2\rho)(3+2\rho)}, \\
\label{eq:573}
\mathbb{COV}[|\mathcal G_i \setminus\mathcal G_j|,|\mathcal G_j
\setminus\mathcal G_k|] & = &
\frac{\theta^2}{(1+\rho)^2(1+2\rho)(3+2\rho)}, \\
\label{eq:574}
\mathbb{COV}[|\mathcal G_i \setminus\mathcal G_j|,|\mathcal G_k
\setminus\mathcal G_j|] & = &
\frac{\theta^2}{(1+\rho)^2(1+2\rho)(3+2\rho)}\nonumber\\[-8pt]\\[-8pt]
&&{} +
\frac{\theta}{(1+\rho)(2+\rho)}.\nonumber
\end{eqnarray}
\end{proposition}
\begin{pf}
We use Proposition \ref{P:main} again. Let $\mathcal T$ be the tree
connecting three individuals $i$, $j$ and $k$, that is, $\mathcal L =
\{i,j,k\}$. Assume the random times $T_2,T_3$ during which the
coalescent has 2, 3 lines, respectively, and one of the three possible
tree topologies, illustrated in Figure \ref{fig:3cases},
%
%
\begin{figure}

\includegraphics{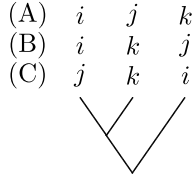}

\caption{The 3 cases for a genealogical tree
connecting three individuals $i,j,k$.}\label{fig:3cases}
\end{figure}
are given. We use
%
%
\begin{eqnarray}
\label{eq:cov1213}
&&\mathbb{COV}[|\mathcal G_i \setminus\mathcal G_j|,|\mathcal G_i
\setminus\mathcal G_k|] \nonumber\\
&&\qquad =
\mathbb{COV} \bigl[\mathbb E[|\mathcal G_i \setminus\mathcal G_j| |
\mathcal T],\mathbb
E[|\mathcal G_i \setminus\mathcal G_k| |\mathcal T] \bigr]
\\
&&\qquad\quad{}
+ \mathbb E \bigl[
\mathbb{COV}[|\mathcal G_i \setminus\mathcal G_j|,|\mathcal G_i
\setminus\mathcal G_k| |\mathcal T] \bigr]\nonumber
\end{eqnarray}
and similar equalities for the other cases. We compute both parts of
the right-hand side separately. For the first part we need to
calculate $\mathbb E [ |\mathcal G_i \setminus\mathcal G_j| | \mathcal
T ]$
depending on~$\mathcal T$:
\begin{enumerate}
\item$\mathcal T \in\{(A)\}$:
\[
\mathbb E [ |\mathcal G_i \setminus\mathcal G_j| | \mathcal T ] =
\int_0^{2T_3} \frac{\theta}{2} e^{-{\rho}/{2}t}\,dt =
\frac{\theta}{\rho}(1-e^{-\rho T_3});
\]
\item$\mathcal T \in\{(B),(C)\}$:
\[
\mathbb E [ |\mathcal G_i \setminus\mathcal G_j| | \mathcal T ] =
\int_0^{2T_2 + 2T_3} \frac{\theta}{2}
e^{-{\rho}/{2}t}\,dt = \frac{\theta}{\rho}\bigl(1-e^{-\rho(T_2 +
T_3)}\bigr).
\]
\end{enumerate}
Replacing the pair $ij$ in the last to expressions by $ik, jk, ki$ or
$kj$ leads to the same possibilities arising in the genealogies $(A),
(B), (C)$. We collect all possibilities in Table
\ref{tab:three}.\vadjust{\goodbreak}

%
%
\begin{table}
\tablewidth=225pt
\caption{The three different tree topologies from Figure
\protect\ref{fig:3cases} give rise to two different terms
for the conditional expectation of a pair of leaves,
depending on the labeling of the pair}\label{tab:three}
\begin{tabular*}{\tablewidth}{@{\extracolsep{\fill}}cccc@{}}
\hline
& $\bolds{(A)}$ & $\bolds{(B)}$ & $\bolds{(C)}$\\
\hline
$ij$ & 1. & 2. & 2.\\
$ji$ & 1. & 2. & 2.\\
$ik$ & 2. & 1. & 2.\\
$jk$ & 2. & 2. & 1.\\
$ki$ & 2. & 1. & 2.\\
$kj$ & 2. & 2. & 1.\\
\hline
\end{tabular*}
\end{table}

In Proposition \ref{P:56} we have seen that $\mathbb E [ \mathbb E
[ |\mathcal G_i \setminus\mathcal G_j| | \mathcal T ] = \mathbb
E[|\mathcal G_i \setminus\mathcal G_j|] =
\frac{\theta}{1+\rho}$ and therefore,
\begin{eqnarray*}
&&
\mathbb{COV} \bigl[\mathbb E[|\mathcal G_i \setminus\mathcal G_j| |
\mathcal T],\mathbb
E[|\mathcal G_i \setminus\mathcal G_k| |\mathcal T] \bigr] \\
&&\qquad= \frac23
\frac{\theta^2}{\rho^2} \mathbb{E} \bigl[ \bigl(1-e^{-\rho(T_2+T_3)}\bigr)
(1-e^{-\rho T_3}) \bigr] \\
&&\qquad\quad{} + \frac13 \frac{\theta^2}{\rho^2}
\mathbb{E} \bigl[\bigl(1-e^{-\rho(T_2+T_3)}\bigr)^2 \bigr] - \frac{\theta^2}{(1+\rho
)^2}\\
&&\qquad= \frac{\theta^2}{3\rho^2} \biggl( \frac{2\rho}{3+\rho} -
\frac{6}{(1+\rho)(3+\rho)} + \frac{6}{(1+\rho)(3+2\rho)} \\
&&\qquad\quad\hspace*{30.3pt}{}
+ 1 -
\frac{6}{(1+\rho)(3+\rho)}+\frac{1}{(1+2\rho)(3+2\rho)} \biggr) -
\frac{\theta^2}{(1+\rho)^2} \\
&&\qquad=
\frac{\theta^2}{(1+\rho)^2(1+2\rho)(3+2\rho)}.
\end{eqnarray*}
Note that this equation also holds for the other three cases in
Proposition \ref{P:57}, that is, we have computed the first term in
(\ref{eq:cov1213}) for all combinations of $i,j,k$ arising in the
proposition.

Let us now consider the second part of (\ref{eq:cov1213}). From
Corollary \ref{Cor:main} we see that
\[
\mathbb E \bigl[ \mathbb{COV}[|\mathcal G_i \setminus\mathcal
G_j|,|\mathcal G_k \setminus\mathcal G_i| |\mathcal T] \bigr]
= \mathbb E \bigl[ \mathbb{COV}[|\mathcal G_i \setminus\mathcal
G_j|,|\mathcal G_j \setminus\mathcal G_k| |\mathcal T] \bigr]
= 0,
\]
which already gives assertions (\ref{eq:572}) and
(\ref{eq:573}). Moreover, Corollary \ref{Cor:main} gives
\begin{eqnarray*}
\mathbb E \bigl[ \mathbb{COV}[|\mathcal G_i \setminus\mathcal
G_j|,|\mathcal G_i \setminus\mathcal G_k| |\mathcal T] \bigr]
& = & \mathbb E \bigl[ \mathbb V [ |\mathcal G_i \setminus
(\mathcal G_j \cup\mathcal G_k)| |\mathcal T ] \bigr],\\
\mathbb E \bigl[ \mathbb{COV}[|\mathcal G_i \setminus\mathcal
G_j|,|\mathcal G_k \setminus\mathcal G_j| |\mathcal T] \bigr]
& = & \mathbb E \bigl[ \mathbb V [ |\mathcal G_i \cap\mathcal G_k \setminus
\mathcal G_j | |\mathcal T ] \bigr].
\end{eqnarray*}
From Proposition \ref{P:main} we know that for given
$\mathcal T, |\mathcal G_i \setminus(\mathcal G_j \cup\mathcal
G_k)|$ and $|\mathcal G_i \cap\mathcal G_k \setminus\mathcal G_j|$
are Poisson distributed. Note that $|\mathcal G_i \setminus(\mathcal
G_j \cup\mathcal G_k)|$ is the number of genes present in $i$, but
not in $j$ and $k$. Recalling that $G_k^{(n)}$ denotes the number of
genes present in $k$ out of $n$ individuals, it is clear that $\mathbb
E[|\mathcal G_i \setminus(\mathcal G_j \cup\mathcal G_k)| ] =
\frac13 \mathbb E[G_1^{(3)}]$, and so using Theorem \ref{T5}
\begin{eqnarray*}
\mathbb E \bigl[ \mathbb V [ |\mathcal G_i \setminus(\mathcal G_j
\cup\mathcal G_k)| |\mathcal T ] \bigr] &=& \mathbb E \bigl[
\mathbb E [ |\mathcal G_i \setminus(\mathcal G_j \cup\mathcal
G_k)| |\mathcal T ] \bigr] \\
&=& \mathbb E[|\mathcal G_i \setminus
(\mathcal G_j \cup\mathcal G_k)| ] = \frac{\theta}{2+\rho}.
\end{eqnarray*}
Equivalently, with $\mathbb E[|\mathcal G_i \cap\mathcal G_k
\setminus\mathcal G_j |] = \frac13 \mathbb E[G_2^{(3)}]$,
\[
\mathbb E \bigl[ \mathbb V [ |\mathcal G_i \cap\mathcal G_k
\setminus\mathcal G_j | |\mathcal T ] \bigr]
= \mathbb E[|\mathcal G_i \cap\mathcal G_k \setminus\mathcal G_j
|] = \frac{\theta}{(2+\rho)(1+\rho)}.
\]
\upqed\end{pf}

\subsection{\texorpdfstring{Gene content for pairs of two
individuals.}{Gene content for pairs of two individuals}}
\label{S:54}
\begin{proposition}[(Gene content for pairs of two individuals)]
\label{P:58}
For $i,j,k,l\in\{1,\ldots,n\}$ pairwise different
\begin{eqnarray*}
&&\mathbb{COV}[|\mathcal G_i \setminus\mathcal G_j|,|\mathcal G_k
\setminus\mathcal G_l|] \\
&&\qquad= \frac{\theta}{(3+\rho)(2+\rho)}
+ \frac{2\theta^2}{(1+\rho)^2(3+\rho)(1+2\rho)(3+2\rho)}.
\end{eqnarray*}
\end{proposition}
\begin{pf}
The proof is similar to the proof of Proposition
\ref{P:57}. Analogously to (\ref{eq:cov1213}) we use
%
%
\begin{eqnarray}
\label{eq:cov1324}
\mathbb{COV} [|\mathcal G_i \setminus\mathcal G_j|,|\mathcal G_k
\setminus\mathcal G_l|]
&=&\mathbb{COV} \bigl[\mathbb E[|\mathcal G_i \setminus\mathcal G_j| |
\mathcal T],
\mathbb E[|\mathcal G_k \setminus\mathcal G_l| |\mathcal T]
\bigr]\nonumber\\[-8pt]\\[-8pt]
&&{} +
\mathbb E \bigl[
\mathbb{COV}[|\mathcal G_i \setminus\mathcal G_j|,|\mathcal G_k
\setminus\mathcal G_l| |\mathcal T] \bigr].\nonumber
\end{eqnarray}
%
As $\mathbb E [ \mathbb E[|\mathcal G_i \setminus\mathcal G_j| |
\mathcal T ] ]
= \mathbb E[|\mathcal G_i \setminus\mathcal G_j|] = \frac{\theta
}{1+\rho}$ we get that
\begin{eqnarray*}
&&
\mathbb{COV} \bigl[ \mathbb{E}[|\mathcal G_i \setminus\mathcal G_j|
|\mathcal T],
\mathbb E[|\mathcal G_k \setminus\mathcal G_l| |\mathcal T] \bigr]
\\
&&\qquad= \mathbb E \bigl[ \mathbb E[ |\mathcal G_i \setminus\mathcal G_j| |
\mathcal T ]
\cdot\mathbb E[ |\mathcal G_k \setminus\mathcal G_l| | \mathcal T ] \bigr]
- \frac{\theta
^2}{(1+\rho)^2}.
\end{eqnarray*}
Therefore, four different cases occur depending on the topology of the
tree seen in Figure \ref{fig:18cases}:

%
%
\begin{figure}

\includegraphics{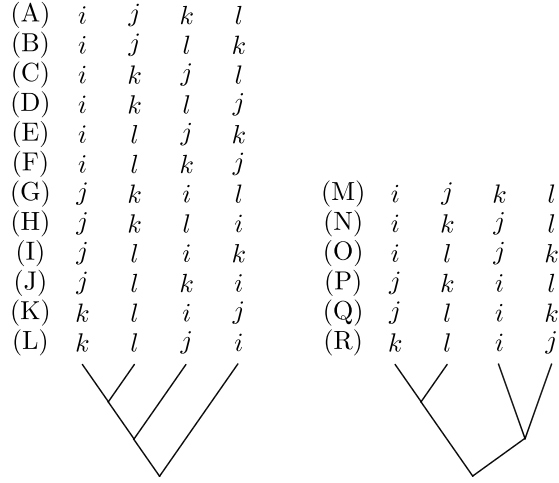}

\caption{The 18 equally probable topologies for
a genealogical tree connecting four individuals $i,j,k,l$.}\label{fig:18cases}
\end{figure}

\begin{enumerate}
\item$\mathcal T \in\{(C),(D),(E),(F),(G),(H),(I),(J)\}$:
\[
\frac{\rho^2}{\theta^2} \mathbb E[ |\mathcal G_i \setminus\mathcal
G_j| | \mathcal T
] \cdot\mathbb E[ |\mathcal G_k \setminus\mathcal G_l| | \mathcal T ]
= \bigl(1 -
e^{-\rho(T_3+T_4)}\bigr)\bigl(1 - e^{-\rho(T_2+T_3+T_4)}\bigr);
\]
\item$\mathcal T \in\{(A),(B),(K),(L)\}$:
\[
\frac{\rho^2}{\theta^2} \mathbb E[ |\mathcal G_i \setminus\mathcal
G_j| | \mathcal T
] \cdot\mathbb E[ |\mathcal G_k \setminus\mathcal G_l| | \mathcal T ]
= (1 -
e^{-\rho T_4})\bigl(1 - e^{\rho(T_2+T_3+T_4)}\bigr);
\]
\item$\mathcal T \in\{(M),(R)\}$:
\[
\frac{\rho^2}{\theta^2} \mathbb E[ |\mathcal G_i \setminus\mathcal
G_j| | \mathcal T
] \cdot\mathbb E[ |\mathcal G_k \setminus\mathcal G_l| | \mathcal T ]
= ( 1 -
e^{-\rho T_4}) \bigl(1 - e^{-\rho(T_3+T_4)}\bigr);
\]
\item$\mathcal T \in\{(N),(O),(P),(Q)\}$:
\[
\frac{\rho^2}{\theta^2} \mathbb E[ |\mathcal G_i \setminus\mathcal
G_j| | \mathcal
T ] \cdot\mathbb E[ |\mathcal G_k \setminus\mathcal G_l| | \mathcal T
] = \bigl(1 -
e^{-\rho(T_2+T_3+T_4)}\bigr)^2.
\]
\end{enumerate}
Hence, with a little help from \textsc{Mathematica},
%
%
\begin{eqnarray}
\label{eq:four1}
&& \mathbb{COV} \bigl[\mathbb E[ |\mathcal G_i \setminus\mathcal G_j| |
\mathcal T ],
\mathbb E[ |\mathcal G_k \setminus\mathcal G_l| | \mathcal T ] \bigr]
\nonumber\\[1pt]
&&\qquad=
\frac{\theta^2}{\rho^2} \biggl(\frac{8}{18} \biggl( 1 -
\frac{18}{(3+\rho)(6+\rho)} -
\frac{18}{(1+\rho)(3+\rho)(6+\rho)} \nonumber\\[1pt]
&&\qquad\quad\hspace*{118.1pt}{}
+
\frac{18}{(1+\rho)(3+2\rho)(6+2\rho)} \biggr)
\nonumber\\[1pt]
&&\qquad\quad\hspace*{19.1pt}{} +
\frac{4}{18} \biggl( 1 - \frac{6}{6+\rho} -
\frac{18}{(1+\rho)(3+\rho)(6+\rho)}\nonumber\\[1pt]
&&\qquad\quad\hspace*{92pt}{} +
\frac{18}{(1+\rho)(3+\rho)(6+2\rho)} \biggr) \nonumber\\[1pt]
&&\qquad\quad\hspace*{19.1pt}{} +
\frac{2}{18} \biggl( 1 - \frac{6}{6+\rho} -
\frac{18}{(3+\rho)(6+\rho)} +\frac{18}{(3+\rho)(6+2\rho)} \biggr)
\end{eqnarray}
\begin{eqnarray}
&&\qquad\quad\hspace*{10pt}{} + \frac{4}{18} \biggl(1 -
\frac{36}{(1+\rho)(3+\rho)(6+\rho)} \hspace*{33.1pt}\hspace*{33.1pt}\nonumber\\
&&\qquad\quad\hspace*{42pt}{}+
\frac{18}{(1+2\rho)(3+2\rho)(6+2\rho)} \biggr) \biggr)-
\frac{\theta^2}{(1+\rho)^2} \hspace*{33.1pt}\hspace*{33.1pt}\nonumber\\
&&\hspace*{11pt}\qquad =
\frac{2\theta^2}{(1+\rho)^2(3+\rho)(1+2\rho)(3+2\rho)}.\hspace*{33.1pt}\hspace*{33.1pt}\nonumber
\end{eqnarray}
%

For the second term, Corollary \ref{Cor:main} gives
\[
\mathbb E \bigl[ \mathbb{COV}[|\mathcal G_i \setminus\mathcal
G_k|,|\mathcal G_j \setminus\mathcal G_l| |\mathcal T] \bigr]
= \mathbb E [ \mathbb V[D_{ij,kl} | \mathcal T ] ]
\]
with $D_{ij,kl}$ as in (\ref{eq:Aijkl}). Given $\mathcal T$, $
D_{ij,kl} $ is Poisson distributed, hence we obtain from
(\ref{eq:EAikl})
%
%
\begin{equation}\label{eq:four2}
\mathbb E [ \mathbb V[D_{ij,kl} | \mathcal T ] ] =
\mathbb E [ \mathbb E[D_{ij,kl} | \mathcal T ] ] =
\frac{\theta}{(3+\rho)(2+\rho)}.
\end{equation}
Combining (\ref{eq:cov1324}) with (\ref{eq:four1}) and
(\ref{eq:four2}) gives the result.
\end{pf}

\subsection{\texorpdfstring{Proof of Theorem \protect\ref{T2}.}{Proof
of Theorem 2}}
\label{S:55}
Using Propositions \ref{P:56}, \ref{P:57} and \ref{P:58}, it is now
easy to prove Theorem \ref{T2}. We obtain
\begin{eqnarray*}
&&
n^2 (n-1)^2 \mathbb V[D] \\
&&\qquad = \sum_{i\neq j} ( \mathbb
V[|\mathcal G_i \setminus\mathcal G_j|] + \mathbb{COV}[|\mathcal
G_i\setminus\mathcal G_j|, |\mathcal G_j\setminus\mathcal
G_i|] )\\
&&\qquad\quad{} + \sum_{i,j,k\ \mathrm{pwd}}
(\mathbb{COV}[|\mathcal G_i\setminus\mathcal G_j|, |\mathcal
G_i\setminus\mathcal G_k|] + \mathbb{COV}[|\mathcal
G_i\setminus\mathcal G_j|, |\mathcal G_j\setminus\mathcal G_k|] \\
&&\qquad\quad\hspace*{49.7pt}{}
+ \mathbb{COV}[|\mathcal
G_i\setminus\mathcal G_j|, |\mathcal G_k\setminus\mathcal G_i|] +
\mathbb{COV}[|\mathcal G_i\setminus\mathcal G_j|, |\mathcal
G_k\setminus\mathcal G_j|] ) \\
&&\qquad\quad{} + \sum_{i,j,k,l\
\mathrm{pwd}} \mathbb{COV}[|\mathcal G_i\setminus\mathcal G_j|,
|\mathcal G_k\setminus\mathcal G_l|] \\
&&\qquad = n(n-1) (\mathbb
V[|\mathcal G_1 \setminus\mathcal G_2|] + \mathbb{COV}[|\mathcal
G_1\setminus\mathcal G_2|, |\mathcal G_2\setminus\mathcal G_1|] )
\\
&&\qquad\quad{} + n(n-1)(n-2)(\mathbb{COV}[|\mathcal
G_1\setminus\mathcal G_2|, |\mathcal G_1\setminus\mathcal G_3|]\\
&&\qquad\quad\hspace*{87pt}{} +
\mathbb{COV}[|\mathcal G_1\setminus\mathcal G_2|, |\mathcal
G_2\setminus\mathcal G_3|] \\
&&\qquad\quad\hspace*{87pt}{} +
\mathbb{COV}[|\mathcal G_1\setminus\mathcal G_2|, |\mathcal
G_3\setminus\mathcal G_1|] \\
&&\qquad\quad\hspace*{88.7pt}{} + \mathbb{COV}[|\mathcal
G_1\setminus\mathcal G_2|, |\mathcal G_3\setminus\mathcal G_2|] )
\\
&&\qquad\quad{} + n(n-1)(n-2)(n-3) \mathbb{COV}[|\mathcal
G_1\setminus\mathcal G_2|, |\mathcal G_3\setminus\mathcal G_4|],
\end{eqnarray*}
and the result follows by some application of \textsc{Mathematica}.

\section{\texorpdfstring{Proof of Theorem \protect\ref{T3}.}{Proof of
Theorem 3}}
We denote by $\mathcal T$ the genealogy connecting the individuals
$i,j,k,l$. As above, we note that $\mathcal T$ is uniquely given by
the random times $T_2, T_3, T_4$ during which the coalescent has
$2,3,4$ lines, respectively, and the tree topology, distinguished by
18 equally probably cases, illustrated in Figure \ref{fig:18cases}. We
use
%
%
\begin{eqnarray}\label{eq:cove}
\mathbb{E}[D_{ij,kl}\cdot D_{ik,jl}] &=& \mathbb
E [\mathbb{COV}[D_{ij,kl}, D_{ik,jl}|\mathcal T] ]
\nonumber\\[-8pt]\\[-8pt]
&&{} +
\mathbb{E} \bigl[\mathbb{E}[D_{ij,kl}|\mathcal T ] \cdot\mathbb
E [D_{ik,jl}|\mathcal T] \bigr]\nonumber
\end{eqnarray}
and note that
\[
\mathbb{COV}[D_{ij,kl}, D_{ik,jl}|\mathcal T] = 0
\]
by Corollary \ref{Cor:main}. So, we are left with computing the
second term in (\ref{eq:cove}). The terms $\mathbb
E[D_{ij,kl}|\mathcal T]$ can take six different values, depending on
$\mathcal T$. We use Proposition \ref{P:main}:

\begin{enumerate}
\item$\mathcal T\in\{(A),(B)\}$:
\begin{eqnarray*}
&& \frac\rho{\theta} \mathbb E[D_{ij,kl}|\mathcal T] \\
&&\qquad=
e^{(-\rho/2) 2T_4} \bigl( 1 - e^{-\rho/2 T_3} +
e^{-\rho/2 T_3}\bigl(1-e^{-\rho/2 (T_3+T_4)}\bigr)\\
&&\qquad\quad\hspace*{114.8pt}{}\times\bigl(1 -
e^{-\rho/2 (2T_2+T_3+T_4)}\bigr) \bigr) \\
&&\qquad= e^{(-\rho/2)
2T_4} - e^{-\rho/2 (2T_3+3T_4)} - e^{-\rho/2
(2T_2+2T_3+3T_4)}\\
&&\qquad\quad{} + e^{-\rho/2(2T_2+3T_3+4T_4)};
\end{eqnarray*}
\item$\mathcal T\in\{(C),(E),(G),(I)\}$:
\begin{eqnarray*}
&& \frac\rho{\theta} \mathbb E[D_{ij,kl}|\mathcal T] =
e^{-\rho/2 (2T_3+2T_4)}(1-e^{-\rho/2T_4})\bigl(1-e^{-\rho/2(2T_2+T_3+T_4)}\bigr) \\
&&\qquad= e^{-\rho/2 (2T_3+2T_4)} - e^{-\rho/2 (2T_3+3T_4)} - e^{-\rho/2 (2T_2+3T_3+3T_4)}\\
&&\qquad\quad{} + e^{-\rho/2 (2T_2+3T_3+4T_4)};
\end{eqnarray*}
\item$\mathcal T\in\{(D),(F),(G),(H),(J), (N), (O), (P),
(Q)\}$:
\begin{eqnarray*}
\frac\rho{\theta} \mathbb E[D_{ij,kl}|\mathcal T]
& = &
e^{-\rho/2 (2T_2+2T_3+2T_4)}(1-e^{-\rho/2T_4})\bigl(1-e^{-\rho/2(T_3+T_4)}\bigr) \\
& = & e^{-\rho/2
(2T_2+2T_3+2T_4)} - e^{-\rho/2 (2T_2+2T_3+3T_4)}\\
&&{} -
e^{-\rho/2 (2T_2+3T_3+3T_4)}
+ e^{-\rho/2 (2T_2+3T_3+4T_4)};
\end{eqnarray*}
\item$\mathcal T\in\{(K),(L)\}$:
\begin{eqnarray*}
\frac\rho{\theta} \mathbb E[D_{ij,kl}|\mathcal T]
& = &
e^{-\rho/2 (2T_2+2T_3+2T_4)}\bigl(1-2e^{-\rho/2(T_3+T_4)} +
e^{-\rho/2(T_3+2T_4)}\bigr) \\
& = & e^{-\rho/2
(2T_2+2T_3+2T_4)} - 2e^{-\rho/2 (2T_2+3T_3+3T_4)}\\
&&{} + e^{-\rho/2 (2T_2+3T_3+4T_4)};
\end{eqnarray*}
\item$\mathcal T\in\{(M)\}$:
\begin{eqnarray*}
\frac\rho{\theta} \mathbb E[D_{ij,kl}|\mathcal T]
& = &
e^{(-\rho/2) 2T_4}\bigl(1 - 2e^{-\rho/2 (2T_2+2T_3+T_4)} +
e^{-\rho/2(2T_2+3T_3+2T_4)}\bigr) \\
& = & e^{(-\rho/2) 2T_4}
- 2e^{-\rho/2 (2T_2+2T_3+3T_4)} + e^{-\rho/2
(2T_2+3T_3+4T_4)};
\end{eqnarray*}
\item$\mathcal T\in\{(R)\}$:
\begin{eqnarray*}
\frac\rho{\theta} \mathbb E[D_{ij,kl}|\mathcal T]
& = &
e^{-\rho/2 (2T_3+2T_4)}\bigl(1 - 2e^{-\rho/2
(2T_2+T_3+T_4)} + e^{-\rho/2(2T_2+T_3+2T_4)}\bigr) \\
& = &
e^{-\rho/2 (2T_3+2T_4)} - 2 e^{-\rho/2 (2T_2 +
3T_3+3T_4)} + e^{-\rho/2 (2T_2+3T_3+4T_4)}.
\end{eqnarray*}
\end{enumerate}

Relabeling $i,j,k,l$ by $i,k,j,l$ changes these cases. Table
\ref{tab:four} gives the responsible terms for $\mathbb
E[D_{ij,kl}|\mathcal T]$ and $\mathbb E[D_{ik,jl}|\mathcal T]$ for all
18 possible tree topologies. We obtain nine cases for which to compute
the second term in (\ref{eq:cove}). We abbreviate $\widetilde e :=
e^{-\rho/2}$.

%
%
\begin{table}
\caption{For every tree topology of Figure \protect\ref{fig:18cases},
the two pairs $ij,kl$ as well as $ik,jl$ fall into one of the six cases
for the conditional expectation; see below
(\protect\ref{eq:cove})}\label{tab:four}
\begin{tabular*}{\tablewidth}{@{\extracolsep{\fill}}cccccccccc@{}}
\hline
& \textbf{(\textit{A})} & \textbf{(\textit{B})} & \textbf{(\textit{C})} & \textbf{(\textit{D})}
& \textbf{(\textit{E})} & \textbf{(\textit{F})} & \textbf{(\textit{G})} & \textbf{(\textit{H})} & \textbf{(\textit{I})} \\
\hline
$ij,kl$ & 1. & 1. & 2. & 3. & 2. & 3. & 2. & 3. & 2. \\
$ik,jl$ & 2. & 3. & 1. & 1. & 3. & 2. & 2. & 3. & 4. \\
\hline
& \textbf{(\textit{J})} & \textbf{(\textit{K})} & \textbf{(\textit{L})} & \textbf{(\textit{M})} & \textbf{(\textit{N})}
& \textbf{(\textit{O})} & \textbf{(\textit{P})} & \textbf{(\textit{Q})} & \textbf{(\textit{R})} \\
\hline
$ij,kl$ & 3. & 4. & 4. & 5. & 3. & 3. & 3. & 3. & 6. \\
$ik,jl$ & 4. & 2. & 3. & 3. & 5. & 3. & 3. & 6. & 3.\\
\hline
\end{tabular*}
\end{table}

$\mathcal T\in\{(A),(C)\}$:
\begin{eqnarray*}
&& \frac{\rho^2}{\theta^2}\mathbb{E}[D_{ij,kl}|\mathcal T] \cdot
\mathbb E[D_{ik,jl}|\mathcal T] \\
&&\qquad= \widetilde e{}^{ 2T_3+4T_4} -
\widetilde e^{ 2T_3+5T_4} - \widetilde e^{ 2T_2+3T_3+5T_4} +
\widetilde e^{ 2T_2+3T_3+6T_4} \\
&&\qquad\quad{} - \widetilde e^{ 4T_3+5T_4} +
\widetilde e^{ 4T_3+6T_4} + \widetilde e^{ 2T_2+5T_3+6T_4} -
\widetilde e^{ 2T_2+5T_3+7T_4} \\
&&\qquad\quad{} - \widetilde e^{
2T_2+4T_3+5T_4} + \widetilde e^{ 2T_2+4T_3+6T_4} + \widetilde
e^{ 4T_2+5T_3+6T_4} - \widetilde e^{ 4T_2+5T_3+7T_4} \\
&&\qquad\quad{} +
\widetilde e^{ 2T_2+5T_3+6T_4} - \widetilde e^{ 2T_2+5T_3+7T_4} -
\widetilde e^{ 4T_2+6T_3+7T_4} + \widetilde e^{ 4T_2+6T_3+8T_4};
\end{eqnarray*}

$\mathcal T\in\{(B),(D)\}$:
\begin{eqnarray*}
&& \frac{\rho^2}{\theta^2} [D_{ij,kl}|\mathcal T] \cdot\mathbb
E[D_{ik,jl}|\mathcal T] \\
&&\qquad= \widetilde e^{ 2T_2+2T_3+4T_4}
- \widetilde e^{ 2T_2+2T_3+5T_4} - \widetilde e^{
2T_2+3T_3+5T_4} + \widetilde e^{ 2T_2+3T_3+6T_4} \\
&&\qquad\quad{} -
\widetilde e^{ 2T_2+4T_3+5T_4} + \widetilde e^{
2T_2+4T_3+6T_4} + \widetilde e^{ 2T_2+5T_3+6T_4} -
\widetilde e^{ 2T_2+5T_3+7T_4} \\
&&\qquad\quad{} - \widetilde e^{
4T_2+4T_3+5T_4} + \widetilde e^{ 4T_2+4T_3+6T_4} +
\widetilde e^{ 4T_2+5T_3+6T_4} - \widetilde e^{
4T_2+5T_3+7T_4} \\
&&\qquad\quad{} + \widetilde e^{ 4T_2+5T_3+6T_4} -
\widetilde e^{ 4T_2+5T_3+7T_4} - \widetilde e^{
4T_2+6T_3+7T_4} + \widetilde e^{ 4T_2+6T_3+8T_4};
\end{eqnarray*}

$\mathcal T\in\{(E),(F)\}$:
\begin{eqnarray*}
\label{eq:coveAC3}
&&
\frac{\rho^2}{\theta^2} [D_{ij,kl}|\mathcal T] \cdot\mathbb
E[D_{ik,jl}|\mathcal T] \\
&&\qquad= \widetilde e^{ 2T_2+4T_3+4T_4}
- \widetilde e^{ 2T_2+4T_3+5T_4} - \widetilde e^{
2T_2+5T_3+5T_4} + \widetilde e^{ 2T_2+5T_3+6T_4} \\
&&\qquad\quad{} -
\widetilde e^{ 2T_2+4T_3+5T_4} + \widetilde e^{
2T_2+4T_3+6T_4} + \widetilde e^{ 2T_2+5T_3+6T_4} -
\widetilde e^{ 2T_2+5T_3+7T_4}\\
&&\qquad\quad{} - \widetilde e^{
4T_2+5T_3+5T_4} + \widetilde e^{ 4T_2+5T_3+6T_4} +
\widetilde e^{ 4T_2+6T_3+6T_4} - \widetilde e^{
4T_2+6T_3+7T_4}\\
&&\qquad\quad{} + \widetilde e^{ 4T_2+5T_3+6T_4} -
\widetilde e^{ 4T_2+5T_3+7T_4} - \widetilde e^{
4T_2+6T_3+7T_4} + \widetilde e^{ 4T_2+6T_3+8T_4};
\end{eqnarray*}

$\mathcal T\in\{(G)\}$:
\begin{eqnarray*}
&&
\frac{\rho^2}{\theta^2} \mathbb{E} [D_{ij,kl}|\mathcal T] \cdot
\mathbb E[D_{ik,jl}|\mathcal T] \\
&&\qquad= \widetilde e^{
4T_3+4T_4} - \widetilde e^{ 4T_3+5T_4} - \widetilde e^{
2T_2+5T_3+5T_4} + \widetilde e^{ 2T_2+5T_3+6T_4} \\
&&\qquad\quad{} -
\widetilde e^{ 4T_3+5T_4} + \widetilde e^{ 4T_3+6T_4} +
\widetilde e^{ 2T_2+5T_3+6T_4} - \widetilde e^{
2T_2+5T_3+7T_4}\\
&&\qquad\quad{} - \widetilde e^{ 2T_2+5T_3+5T_4} +
\widetilde e^{ 2T_2+5T_3+6T_4} + \widetilde e^{
4T_2+6T_3+6T_4} - \widetilde e^{ 4T_2+6T_3+7T_4}\\
&&\qquad\quad{} +
\widetilde e^{ 2T_2+5T_3+6T_4} - \widetilde e^{
2T_2+5T_3+7T_4} - \widetilde e^{ 4T_2+6T_3+7T_4} +
\widetilde e^{ 4T_2+6T_3+8T_4};
\end{eqnarray*}

$\mathcal T\in\{(H),(O),(P)\}$:
\begin{eqnarray*}
&&\frac{\rho^2}{\theta^2} \mathbb{E} [D_{ij,kl}|\mathcal T] \cdot
\mathbb E[D_{ik,jl}|\mathcal T] \\
&&\qquad= \widetilde e^{
4T_2+4T_3+4T_4} - \widetilde e^{ 4T_2+4T_3+5T_4} -
\widetilde e^{ 4T_2+5T_3+5T_4} + \widetilde e^{
4T_2+5T_3+6T_4} \\
&&\qquad\quad{} - \widetilde e^{ 4T_2+4T_3+5T_4} +
\widetilde e^{ 4T_2+4T_3+6T_4} + \widetilde e^{
4T_2+5T_3+6T_4} - \widetilde e^{ 4T_2+5T_3+7T_4}\\
&&\qquad\quad{} -
\widetilde e^{ 4T_2+5T_3+5T_4} + \widetilde e^{
4T_2+5T_3+6T_4} + \widetilde e^{ 4T_2+6T_3+6T_4} -
\widetilde e^{ 4T_2+6T_3+7T_4}\\
&&\qquad\quad{} + \widetilde e^{
4T_2+5T_3+6T_4} - \widetilde e^{ 4T_2+5T_3+7T_4} -
\widetilde e^{ 4T_2+6T_3+7T_4} + \widetilde e^{
4T_2+6T_3+8T_4};
\end{eqnarray*}

$\mathcal T\in\{(I),(K)\}$:
\begin{eqnarray*}
&&
\frac{\rho^2}{\theta^2} \mathbb{E} [D_{ij,kl}|\mathcal T] \cdot
\mathbb E[D_{ik,jl}|\mathcal T] \\
&&\qquad= \widetilde e^{
2T_2+4T_3+4T_4} - 2\widetilde e^{ 2T_2+5T_3+5T_4}+
\widetilde e^{ 2T_2+5T_3+6T_4} \\
&&\qquad\quad{} - \widetilde e^{
2T_2+4T_3+5T_4} + 2\widetilde e^{ 2T_2+5T_3+6T_4}-
\widetilde e^{ 2T_2+5T_3+7T_4} \\
&&\qquad\quad{} - \widetilde e^{
4T_2+5T_3+5T_4} + 2\widetilde e^{ 4T_2+6T_3+6T_4} -
\widetilde e^{ 4T_2+6T_3+7T_4} \\
&&\qquad\quad{} + \widetilde e^{
4T_2+5T_3+6T_4} - 2\widetilde e^{ 4T_2+6T_3+7T_4}+
\widetilde e^{ 4T_2+6T_3+8T_4};
\end{eqnarray*}

$\mathcal T\in\{(J),(L)\}$:
\begin{eqnarray*}
&&
\frac{\rho^2}{\theta^2} \mathbb{E} [D_{ij,kl}|\mathcal T] \cdot
\mathbb E[D_{ik,jl}|\mathcal T] \\
&&\qquad= \widetilde e^{
4T_2+4T_3+4T_4} - 2\widetilde e^{ 4T_2+5T_3+5T_4}+
\widetilde e^{ 4T_2+5T_3+6T_4} \\
&&\qquad\quad{} - \widetilde e^{
4T_2+4T_3+5T_4} + 2\widetilde e^{ 4T_2+5T_3+6T_4} -
\widetilde e^{ 4T_2+5T_3+7T_4} \\
&&\qquad\quad{} - \widetilde e^{
4T_2+5T_3+5T_4} + 2\widetilde e^{ 4T_2+6T_3+6T_4}-
\widetilde e^{ 4T_2+6T_3+7T_4} \\
&&\qquad\quad{} + \widetilde e^{
4T_2+5T_3+6T_4} - 2\widetilde e^{ 4T_2+6T_3+7T_4}+
\widetilde e^{ 4T_2+6T_3+8T_4};
\end{eqnarray*}

$\mathcal T\in\{(M),(N)\}$:
\begin{eqnarray*}
&&\frac{\rho^2}{\theta^2} \mathbb{E} [D_{ij,kl}|\mathcal T] \cdot
\mathbb E[D_{ik,jl}|\mathcal T] \\
&&\qquad= \widetilde e^{
2T_2+2T_3+4T_4} - 2\widetilde e^{ 4T_2+4T_3+5T_4} +
\widetilde e^{ 4T_2+5T_3+6T_4} \\
&&\qquad\quad{} - \widetilde e^{
2T_2+2T_3+5T_4} + 2\widetilde e^{ 4T_2+4T_3+6T_4} -
\widetilde e^{ 4T_2+5T_3+7T_4}\\
&&\qquad\quad{} - \widetilde e^{
2T_2+3T_3+5T_4} + 2\widetilde e^{ 4T_2+5T_3+6T_4} -
\widetilde e^{ 4T_2+6T_3+7T_4}\\
&&\qquad\quad{} + \widetilde e^{
2T_2+3T_3+6T_4} - 2\widetilde e^{ 4T_2+5T_3+7T_4} +
\widetilde e^{ 4T_2+6T_3+8T_4};
\end{eqnarray*}

$\mathcal T\in\{(Q),(R)\}$:
\begin{eqnarray*}
&&
\frac{\rho^2}{\theta^2} \mathbb{E} [D_{ij,kl}|\mathcal T] \cdot
\mathbb E[D_{ik,jl}|\mathcal T] \\
&&\qquad= \widetilde e^{
2T_2+4T_3+4T_4} - 2 \widetilde e^{ 4T_2 + 5T_3+5T_4} +
\widetilde e^{ 4T_2+5T_3+6T_4}\\
&&\qquad\quad{} - \widetilde e^{
2T_2+4T_3+5T_4} + 2 \widetilde e^{ 4T_2 + 5T_3+6T_4} -
\widetilde e^{ 4T_2+5T_3+7T_4}\\
&&\qquad\quad{} - \widetilde e^{
2T_2+5T_3+5T_4} + 2 \widetilde e^{ 4T_2 + 6T_3+6T_4} -
\widetilde e^{ 4T_2+6T_3+7T_4}\\
&&\qquad\quad{} + \widetilde e^{
2T_2+5T_3+6T_4} - 2 \widetilde e^{ 4T_2 + 6T_3+7T_4} +
\widetilde e^{ 4T_2+6T_3+8T_4}.
\end{eqnarray*}
Combining the last equations and using \textsc{Mathematica}, we obtain
the desired result.




\section{\texorpdfstring{Proof of Theorems \protect\ref{T4} and
\protect\ref{T5}.}{Proof of Theorems 4 and 5}}

Theorems \ref{T4} and \ref{T5} provide expectations for the size of
the dispensable genome and the gene frequency spectrum,
respectively. Recalling that $G_k^{(n)}$ is the number of genes in
frequency $k=1,\ldots,n$ in a sample of size $n$, it is clear that
\[
G = \sum_{k=1}^n G_k^{(n)},
\]
where $G$ is given by (\ref{eq:G}). In addition, if $\mathbb
E[G_1^{(i)}]$, $i=1,\ldots,$ is known from Theorem \ref{T5},
\[
\mathbb E[G] = \sum_{i=1}^n \frac1i
\mathbb E\bigl[G_1^{(i)}\bigr] = \sum_{i=1}^n \frac{1}{i} \frac{\theta
i}{i-1+\rho}
= \theta\sum_{i=0}^{n-1} \frac{1}{i+\rho}.
\]
Hence, the result for the expected number of genes in Theorem \ref{T4}
can easily be proved once we have established Theorem
\ref{T5}. However, we take an alternative route and give an
independent proof of Theorem \ref{T4}.

\subsection{\texorpdfstring{An independent proof of Theorem
\protect\ref{T4}.}{An independent proof of Theorem 4}}

Recall the survival function $p_{\mathcal T}$ from Definition
\ref{def:tree}. Consider the coalescent, denoted by $\mathcal T$,
started with $n$ lines, rooted at $r$, the most recent common ancestor
of the sample. As shown in Proposition \ref{P:52}, $p_{\mathcal T}\dvtx
\mathcal T\to[0,1]$ gives the probability that a mutation that arises
at $t\in\mathcal T$ is not lost in at least one leaf. Hence, given
$\mathcal T$, we find that $G \sim\mathcal{POI}(\pi(\mathcal T))$
with
\[
\pi(\mathcal T) := \frac\theta2\int_{\mathcal T} p_{\mathcal
T}(t) \,dt.
\]
Next, we consider a random coalescent $\mathcal T$ with
additional \textit{loss events} at rate $\frac\rho2$ along the
tree. We say that $t\in\mathcal T$ is unlost if there is a leaf in
$i\in\mathcal T$ such that the path $[t,i]$ is not hit by a loss
event. Given $\mathcal T$, note that
\[
p_{\mathcal T}(t) = \mathbb P[t \mbox{ unlost}|\mathcal T]
\]
by Proposition \ref{P:52}.

To prove (\ref{eq:T3a}), we write immediately, using the above
arguments,
%
%
\begin{eqnarray}
\label{eq:jkl}
\mathbb E[G] & = & \mathbb E \biggl[\frac\theta2 \int_{\mathcal T}
p_{\mathcal T}(t) \,dt \biggr] = \mathbb E \biggl[\frac\theta2 \int
_{\mathcal T} 1(\mbox{$t$
not lost}) \,dt \biggr] \nonumber\\[-8pt]\\[-8pt]
& = & \frac\theta2 \mathbb E [\mbox{length of unlost lines in
}\mathcal T].\nonumber
\end{eqnarray}

To compute the expected length of unlost lines we note that all lines
are unlost near the leaves. The number of unlost lines decreases
either by a coalescence or by a gene loss event. When there are $k$
unlost lines left, any line is lost at rate $\frac\rho2$ and two
lines coalesce at rate ${k\choose2}$. Hence, the time until there are
$k-1$ unlost lines is $\exp( \frac\rho2 k + {k\choose2} )$
distributed. Thus,
\[
\mathbb E[G] = \frac\theta2 \sum_{k=1}^n
\frac{k}{{k\choose2} + \rho/2 k} = \theta\sum_{k=1}^n
\frac{1}{k-1+\rho} = \theta\sum_{k=0}^{n-1} \frac{1}{k+\rho},
\]
and we have proven (\ref{eq:T3a}).

Next we show how to obtain the recursion for $\mathbb V[G]$ given in
(\ref{eq:T3b}). Using the fact that, given $\mathcal T$, the number of
genes $G$ is Poisson distributed with rate $\pi(\mathcal T)$,
\begin{eqnarray*}
\mathbb V[G] & = & \mathbb E [\mathbb E[G^2|\mathcal T] ] -
\mathbb E[G]^2 = \mathbb E[\pi(\mathcal T) + \pi(\mathcal T)^2] -
\mathbb E[G]^2 \\
& = &\mathbb E[G] - \mathbb E[G]^2 +
\frac{\theta^2}4 \mathbb E \biggl[ \int\int p_{\mathcal T}(s)
p_{\mathcal T}(t) \,ds \,dt \biggr].
\end{eqnarray*}
Since $\mathbb E[G]$ is known, it remains to compute the last term in
the last display. Consider two independent Poisson processes $\mathcal
P_1$ and $\mathcal P_2$ along the tree $\mathcal T$, each at
rate~$\frac\rho2$, describing gene loss. As above, we say that a point
$s\in\mathcal T$ is $k$-unlost if there is a leaf $i\in\mathcal T$
such that the path $[s,i]$ is not hit by an event in $\mathcal
P_k$. We denote by $L_k$ the length of $k$-unlost points in $\mathcal
T$, $k=1,2$. Using the same reasoning as in~(\ref{eq:jkl}),
\[
\mathbb E \biggl[\int\int p_{\mathcal T}(s) p_{\mathcal T}(t) \,ds \,dt \biggr]
= \mathbb E[L_1L_2].
\]
The latter expectation can be derived via the
following construction: in the tree $\mathcal T$ with the two
independent Poisson loss processes $\mathcal P_1$ and $\mathcal P_2$,
denote by $K_1(\tau)$ the number of lines which are both $1$- and
$2$-unlost some distance $\tau$ from the treetop, by $K_2(\tau)$ the
number of lines which are $2$-lost but $1$-unlost and by $K_3(\tau)$
the number of lines which are $1$-lost but $2$-unlost by time
$\tau$. Clearly,
\[
L_1 = \int_0^\infty\bigl(K_1(\tau)+K_2(\tau)\bigr)\,d\tau,\qquad
L_2 = \int_0^\infty\bigl(K_1(\tau)+K_3(\tau)\bigr)\,d\tau.
\]
In addition,
$K=(K(\tau))_{\tau\geq0}= (K_1(\tau),K_2(\tau),K_3(\tau))_{\tau
\geq
0}$ is a Markov jump process with the following rates from
$(k_1,k_2,k_3)$ to

\vspace*{5pt}

\begin{center}
\begin{tabular}{@{}ll@{}}
\hline
New state & At rate \\
\hline
$\underline k_1' = (k_1-1,k_2,k_3)$ &
$\lambda_1={k_1\choose2}$ \\[2pt]
$\underline k_2' = (k_1,k_2-1,k_3)$ &
$\lambda_2={k_2\choose2} + k_1k_2 + \frac\rho2 k_2$ \\[2pt]
$\underline k_3' = (k_1,k_2,k_3-1)$ &
$\lambda_3={k_3\choose2}
+ k_1k_3 + \frac\rho2 k_3$ \\[2pt]
$\underline k_4' = (k_1+1,k_2-1,k_3-1)$ &
$\lambda_4=k_2k_3$\\[2pt]
$\underline k_5' = (k_1-1,k_2+1,k_3)$ &
$\lambda_5=\frac\rho2 k_1$ \\[2pt]
$\underline k_6' = (k_1-1,k_2,k_3+1)$ &
$\lambda_6=\frac\rho2 k_1$\\
\hline
\end{tabular}
\end{center}

\mbox{}


Note that $\underline k_1',\ldots,\underline k_6',
\lambda_1,\ldots,\lambda_6$ are as defined in Definition \ref{def:gh}.
In each transition, $2k_1+k_2+k_3$ is not increasing and therefore
hits 0 after a finite number of transitions.

To obtain (\ref{eq:T3b}), define the process $K$ with $K(0) =
\underline k := (k_1, k_2, k_3)$ and define
\[
L_1^{\underline k} =
\int_0^\infty\bigl(K_1(\tau)+K_2(\tau)\bigr)\,d\tau,\qquad
L_2^{\underline k} = \int_0^\infty\bigl(K_1(\tau)+K_3(\tau)\bigr)\,d\tau
\]
such that $L_i
\stackrel{d}{=} L_i^{(n,0,0)}$, $i=1,2$. We claim that
\[
g_{\underline k}:= \mathbb E[L_1^{\underline k}L_2^{\underline k}]
\]
satisfies (\ref{eq:defgh1})
as well as the recursion (\ref{eq:defgh2}).

First, given $k_1+k_2=1$ ($k_1+k_3=1$), there is only one $1$-unlost
line ($2$-unlost line) in $\mathcal T$. In this case, this line can
only be lost by an event in $\mathcal P_1$ ($\mathcal P_2$),
independent of all other coalescence events. Hence, in this case,
$L_1^{\underline k}$ and $L_2^{\underline k}$ are independent,
$\mathbb E[L_1^{\underline k}]=\frac2\rho$ ($\mathbb
E[L_2^{\underline k}]=\frac2\rho$) and $\mathbb E[L_2^{\underline
k}] = h_{k_1+k_3}$ ($\mathbb E[L_1^{\underline k}] =
h_{k_1+k_2}$). Combining these results gives (\ref{eq:defgh1}).

Second we show that $\mathbb E[L_1^{\underline k}L_2^{\underline k}]$
satisfies (\ref{eq:defgh2}). Since $K$ is a jump process, the first
event occurs after an exponential time $T$ with rate
$\overline\lambda$, which is independent of the new state after the
first jump. Conditioning on the first event happening at time~$T$,
\[
L^{\underline k}_1L^{\underline k}_2 = \sum_{i=1}^6 1_{\{\mathrm{new}
\ \mathrm{state}\ \mathrm{is}\ \underline k_i'\}} \bigl((k_1+k_2)T+L^{\underline
k_i'}_1\bigr)\bigl((k_1+k_3)T+L{}^{\underline k_i'}_2\bigr).
\]
Taking expectations on both sides shows that $\mathbb E[L^{\underline
k}_1L^{\underline k}_2]$ satisfies (\ref{eq:defgh2}). This completes
the proof.


\subsection{\texorpdfstring{Proof of Theorem \protect\ref{T5}.}{Proof
of Theorem 5}}

There are several ways to prove Theorem \ref{T5}. We present here two
approaches, one based on diffusion theory, the other one using an urn
model.
\begin{pf*}{Proof of Theorem \ref{T5} based on diffusion theory}
Assume that a gene is present at frequency $X_0$ at time $0$. Then,
$(X_t)_{t\geq0}$ follows the SDE
\[
dX = -\frac{\rho}{2} X \,dt + \sqrt{X(1-X)}\,dW.
\]
Frequency spectra for such diffusions have been obtained by
\citet{Kimura1964}. We follow the arguments given in
\citet{Durrett2008}, Theorem 7.20. Assume we introduce new genes at
frequency $0<\delta<1$ into the population at rate
\[
\frac\theta2 \frac{1}{\phi(\delta)},
\]
where $\mu(x) := -\frac\rho2 x, \sigma^2(x) := x(1-x)$,
\begin{eqnarray*}
\psi(y)  :\!&=& \exp\biggl(-2 \int_0^y
\frac{\mu(z)}{\sigma^2(z)}\,dz \biggr) = \exp\biggl( \rho\int_0^y
\frac{1}{1-z} \,dz \biggr) \\
& = &\exp\bigl(-\rho
\log(1-y)\bigr) = (1-y)^{-\rho},\\
\phi(x) :\!&=& \int_0^x \psi(y) \,dy =
\frac{1}{1-\rho} \bigl(1-(1-x)^{1-\rho} \bigr).
\end{eqnarray*}
This rate is consistent in $\delta$: the number of genes at level
$\varepsilon>\delta$ is $ \frac\theta2
\frac{1}{\phi(\delta)}\frac{\phi(\delta)}{\phi(\varepsilon)}$ since
$\frac{\phi(\delta)}{\phi(\varepsilon)}$ is the probability that the
gene reaches frequency $\varepsilon$ before dying out. Moreover, the
Green function for the diffusion---measuring the time until
eventual loss of the gene---is given by $2\phi(\delta)m(y)$ for
$y>\delta$, where
\[
m(y) = \frac{1}{\sigma^2(y) \psi(y)} = \frac{1}{y(1-y)^{1-\rho}}
\]
is the density of the speed measure of the diffusion. Hence, we find
that the number of genes in frequency $x$ is Poisson with mean
\[
g(x)\,dx := \theta\frac{1}{x(1-x)^{1-\rho}}\,dx.
\]
Now, the theorem follows since
\begin{eqnarray*}
\mathbb E[G_k] & = & \pmatrix{n\cr k} \int_0^1 g(x) x^k(1-x)^{n-k}\,dx =
\theta\pmatrix{n\cr k} \int_0^1 x^{k-1}(1-x)^{n-k-1+\rho}\,dx \\
& = &
\frac\theta k \frac{n\cdots(n-k+1)}{(n-1+\rho)\cdots
(n-k+\rho)}.
\end{eqnarray*}
\upqed\end{pf*}
\begin{pf*}{Proof of Theorem \ref{T5} based on an urn model}
Let $\mathcal T$ be the Kingman coalescent and $\Gamma_{\mathcal T}$
be the tree-indexed Markov chain from Definition
\ref{def:21}. First, we focus on loss events for $du\subseteq I$
along the random tree $\mathcal T$. (We use the infinitesimal symbol
$du$ for notational convenience.) Since $du$ is small, we may safely
assume that there is at most one gene in $du$ present in
$\bigcup_{i=1}^n \mathcal G_i$. Gene loss events in $du$ occur at
constant rate $\frac\rho2$ along each branch. Consider the tree
$\mathcal T$ from the leaves to the root. Lines coalesce with pair
coalescence rate $1$, and any line hits a loss event in $du$ at rate
$\frac\rho2$. Upon a loss event we kill the line off the
tree. The resulting forest is well known from the family
decomposition in the infinite alleles model [e.g.,
\citet{Durrett2008}, page 14]. Using Hoppe's urn, we can also
generate the forest forward in time: consider an urn with one
colored and one black ball. Choose the colored ball with
probability proportional to 1 and the black one with probability
proportional to $\rho$. When choosing a colored ball, put the
chosen ball plus one ball of the same color into the urn. When
choosing the black ball, put the black ball back together with a
ball of a new color. In the next step, again choose any colored
ball with probability proportional to 1 and the black balls with
probability proportional to $\rho$. Proceed until there are $n$
colored balls in the urn. Note that, given there are $i$ colored
balls in the urn, the chance that the next chosen ball is colored
is $\frac{i}{i+\rho} = \frac{{i+1\choose2}}{{i+1\choose2} +
(i+1)\rho/2}$, that is, the chance equals the probability that
two among $i+1$ lines coalesce and are not killed off the tree by a
gene loss event.

To obtain the correct branch lengths in the tree, when there are $i$
colored balls in the urn, wait an exponential time with rate $\frac
i2 (i-1+\rho)$ until adding the next colored ball. This waiting time
equals the time the coalescent stays with $i$ lines, when pairs
coalesce at rate 1 and single lines are killed at rate $\frac\rho
2$. Hence, by this procedure, balls with the same color belong to
the same tree in the forest, and the time the forest spends with $i$
lines is the same as viewing the coalescent backward in time.

So far, Hoppe's urn only described gene loss of the single gene
$u$. Let us add gene gain of a gene in $du$ to the
description. During the evolution of Hoppe's urn, which comes with
its exponential waiting times, mark all colored balls at rate
$\frac\theta2 \,du$. When a marked colored ball is chosen, the added
ball is again marked. Here, a mark stands for the presence of the
considered gene along the corresponding ancestral line. Since $du$
is small, there is at most one mark along the forest.

For the forest given by the marked Hoppe's urn, we distinguish the
times $T_1, \ldots, T_n$ when there are $1,\ldots,n$ lines present. We say
that line $l$ during $T_i$ is of size $k$ iff the ball belonging to
this line produces exactly $k-1$ offspring until the urn
finishes. Hence,
\begin{eqnarray*}
\mathbb E[G_k] & = &\int_I \mathbb E[du\in\mathcal G_i \mbox{ for
exactly $k$ different $i$}] \\
& = &\sum_{i=1}^n \sum_{l=1}^i
\mathbb P[\mbox{$l$th line during $T_i$ is of size $k$}] \\
&&\hspace*{28.1pt}{} \times\int\mathbb P[\mbox{mark in $du$ on $l$th
line during $T_i$}]
\end{eqnarray*}
and
\[
\mathbb P[\mbox{mark in $du$ on $l$th line during $T_i$}] =
\frac{\theta/2 \,du}{ i/2(i-1+\rho)} =
\frac{\theta}{i(i-1+\rho)}\,du.
\]
Let us turn to the probability that the $l$th line during $T_i$ is
of size $k$. The reasoning below is well known from P\`{o}lya urn
models. When starting with $i-1$ unmarked and one marked lines,
there are ${n-i\choose k-1}$ possibilities at what times $k-1$
marked balls are added when $n-i$ balls are added to the urn in
total. For any of these possibilities, the probability is
\[
\frac{(k-1)!(i-1+\rho)\cdots(n-k-1+\rho)}{(i+\rho)\cdots(n-1+\rho)}.
\]
Putting everything together,
\begin{eqnarray*}
\mathbb E[G_k] & = &\sum_{i=1}^n i \pmatrix{n-i\cr k-1}
\frac{(k-1)!(i-1+\rho)\cdots(n-k-1+\rho)}{(i+\rho)\cdots(n-1+\rho)}
\frac{\theta}{i(i-1+\rho)} \\
& = &
\frac{\theta}{k}\frac{k!}{(n-k+\rho)\cdots(n-1+\rho)}
\underbrace{\sum_{i=1}^n \pmatrix{n-i\cr k-1}}_{= {n\choose k} } \\
& = &
\frac{\theta}{k}\frac{n\cdots(n-k+1)}{(n-k+\rho)\cdots(n-1+\rho)},
\end{eqnarray*}
and we are done.
\end{pf*}

\section*{\texorpdfstring{Acknowledgments.}{Acknowledgments}}
We thank Daniel Huson for pointing out the reference
[\citet{Huson2004Bioinformatics15044248}] and we are grateful to
Cornelia Borck, Andrej Depperschmidt and Bernhard Haubold for comments
on our manuscript.


%
\printaddresses


\begin{thebibliography}{30}

\bibitem[\protect\citeauthoryear{Bentley}{2009}]{Bentley2009}
\begin{barticle}[auto:SpringerTagBib|2009-01-14|16:51:27]
\bauthor{\bsnm{Bentley},~\bfnm{S.}\binits{S.}}
(\byear{2009}).
\btitle{Sequencing the species pan-genome}.
\bjournal{Nature Rev. Microbiol.}
\bvolume{7}
\bpages{258--259}.
\end{barticle}
\endbibitem

\bibitem[\protect\citeauthoryear{Dufresne et al.}{2008}]{Dufresne2008GenomeBiol18507822}
\begin{barticle}[auto:SpringerTagBib|2009-01-14|16:51:27]
\bauthor{\bsnm{Dufresne},~\bfnm{A.}\binits{A.}},
  \bauthor{\bsnm{Ostrowski},~\bfnm{M.}\binits{M.}},
  \bauthor{\bsnm{Scanlan},~\bfnm{D.~J.}\binits{D.~J.}},
  \bauthor{\bsnm{Garczarek},~\bfnm{L.}\binits{L.}},
  \bauthor{\bsnm{Mazard},~\bfnm{S.}\binits{S.}},
  \bauthor{\bsnm{Palenik},~\bfnm{B.~P.}\binits{B.~P.}},
  \bauthor{\bsnm{Paulsen},~\bfnm{I.~T.}\binits{I.~T.}},
  \bauthor{\bparticle{de~}\bsnm{Marsac},~\bfnm{N.~T.}\binits{N.~T.}},
  \bauthor{\bsnm{Wincker},~\bfnm{P.}\binits{P.}},
  \bauthor{\bsnm{Dossat},~\bfnm{C.}\binits{C.}},
  \bauthor{\bsnm{Ferriera},~\bfnm{S.}\binits{S.}},
  \bauthor{\bsnm{Johnson},~\bfnm{J.}\binits{J.}},
  \bauthor{\bsnm{Post},~\bfnm{A.~F.}\binits{A.~F.}},
  \bauthor{\bsnm{Hess},~\bfnm{W.~R.}\binits{W.~R.}} \AND
  \bauthor{\bsnm{Partensky},~\bfnm{F.}\binits{F.}}
  (\byear{2008}).
\btitle{Unraveling the genomic mosaic of a ubiquitous genus of
  marine cyanobacteria}.
  \bjournal{Genome Biol.}
  \bvolume{9}
  \bpages{R90}.
\end{barticle}
\endbibitem

\bibitem[\protect\citeauthoryear{Durrett}{2008}]{Durrett2008}
\begin{bbook}[mr]
\bauthor{\bsnm{Durrett},~\bfnm{Richard}\binits{R.}}
(\byear{2008}).
\btitle{Probability Models for {DNA} Sequence Evolution},
\bedition{2nd} ed.
\bpublisher{Springer}, \baddress{New York}.
\bid{mr={2439767}}
\end{bbook}
\endbibitem

\bibitem[\protect\citeauthoryear{Durrett and Popovic}{2009}]{DurrettPopovic2007}
\begin{barticle}[mr]
\bauthor{\bsnm{Durrett},~\bfnm{Rick}\binits{R.}} \AND
  \bauthor{\bsnm{Popovic},~\bfnm{Lea}\binits{L.}}
(\byear{2009}).
\btitle{Degenerate diffusions arising from gene duplication models}.
\bjournal{Ann. Appl. Probab.}
\bvolume{19}
\bpages{15--48}.
\bid{doi={10.1214/08-AAP530}, mr={2498670}}
\end{barticle}
\endbibitem

\bibitem[\protect\citeauthoryear{Dykhuizen and Green}{1991}]{DykhuizenGreen1991}
\begin{barticle}[auto:SpringerTagBib|2009-01-14|16:51:27]
\bauthor{\bsnm{Dykhuizen},~\bfnm{D.~E.}\binits{D.~E.}} \AND
  \bauthor{\bsnm{Green},~\bfnm{L.}\binits{L.}}
  (\byear{1991}).
\btitle{{R}ecombination in \textit{Escherichia coli} and the
  definition of biological species}.
  \bjournal{J. Bacteriol.}
  \bvolume{173}
\bpages{7257--7268}.
\end{barticle}
\endbibitem

\bibitem[\protect\citeauthoryear{Ehrlich et al.}{2005}]{Ehrlich2005}
\begin{barticle}[auto:SpringerTagBib|2009-01-14|16:51:27]
\bauthor{\bsnm{Ehrlich},~\bfnm{G.~D.}\binits{G.~D.}},
  \bauthor{\bsnm{Hu},~\bfnm{F.~Z.}\binits{F.~Z.}},
  \bauthor{\bsnm{Shen},~\bfnm{K.}\binits{K.}},
  \bauthor{\bsnm{Stoodley},~\bfnm{P.}\binits{P.}} \AND
  \bauthor{\bsnm{Post},~\bfnm{J.~C.}\binits{J.~C.}}
  (\byear{2005}).
\btitle{{B}acterial plurality as a general mechanism driving
  persistence in chronic infections}.
  \bjournal{Clin. Orthop. Relat. Res.}
  \bvolume{437}
  \bpages{20--24}.
\end{barticle}
\endbibitem

\bibitem[\protect\citeauthoryear{Evans, Shvets and
  Slatkin}{2007}]{EvansShvetsSlatkin2007}
\begin{barticle}[auto:SpringerTagBib|2009-01-14|16:51:27]
\bauthor{\bsnm{Evans},~\bfnm{S.}\binits{S.}},
  \bauthor{\bsnm{Shvets},~\bfnm{S.}\binits{S.}} \AND
  \bauthor{\bsnm{Slatkin},~\bfnm{M.}\binits{M.}}
  (\byear{2007}).
\btitle{Non-equlibrium theory of the allele frequency spectrum}.
  \bjournal{Theo. Pop. Biol.}
  \bvolume{71}
  \bpages{109--119}.
\end{barticle}
\endbibitem

\bibitem[\protect\citeauthoryear{Ewens}{2004}]{Ewens2004}
\begin{bbook}[mr]
\bauthor{\bsnm{Ewens},~\bfnm{Warren~J.}\binits{W.~J.}}
(\byear{2004}).
\btitle{Mathematical Population Genetics. {I}. Theoretical Introduction},
\bedition{2nd} ed.
\bseries{Interdisciplinary Applied Mathematics}
\bvolume{27}.
\bpublisher{Springer}, \baddress{New York}.
\bid{mr={2026891}}
\end{bbook}
\endbibitem

\bibitem[\protect\citeauthoryear{Fraser, Hanage and
  Spratt}{2007}]{pmid17255503}
\begin{barticle}[pbm]
\bauthor{\bsnm{Fraser},~\bfnm{Christophe}\binits{C.}},
  \bauthor{\bsnm{Hanage},~\bfnm{William~P.}\binits{W.~P.}} \AND
  \bauthor{\bsnm{Spratt},~\bfnm{Brian~G.}\binits{B.~G.}}
(\byear{2007}).
\btitle{Recombination and the nature of bacterial speciation}.
\bjournal{Science}
\bvolume{315}
\bpages{476--480}.
\bid{pii={315/5811/476}, doi={10.1126/science.1127573}, pmid={17255503},
  pmcid={PMC2220085}, mid={UKMS1405}}
\end{barticle}
\endbibitem

\bibitem[\protect\citeauthoryear{Fu}{1995}]{Fu1995}
\begin{barticle}[auto:SpringerTagBib|2009-01-14|16:51:27]
\bauthor{\bsnm{Fu},~\bfnm{Y.~X.}\binits{Y.~X.}}
(\byear{1995}).
\btitle{{S}tatistical properties of segregating sites}.
\bjournal{Theo. Pop. Biol.}
\bvolume{48}
\bpages{172--197}.
\end{barticle}
\endbibitem

\bibitem[\protect\citeauthoryear{Griffiths}{2003}]{Griffiths2003}
\begin{barticle}[auto:SpringerTagBib|2009-01-14|16:51:27]
\bauthor{\bsnm{Griffiths},~\bfnm{R.~C.}\binits{R.~C.}}
(\byear{2003}).
\btitle{The frequency spectrum of a mutation and its age, in a
  general diffusion model}.
\bjournal{Theo. Pop. Biol.}
\bvolume{64}
\bpages{241--251}.
\end{barticle}
\endbibitem

\bibitem[\protect\citeauthoryear{Hiller et al.}{2007}]{pmid17675389}
\begin{barticle}[auto:SpringerTagBib|2009-01-14|16:51:27]
\bauthor{\bsnm{Hiller},~\bfnm{N.~L.}\binits{N.~L.}},
  \bauthor{\bsnm{Janto},~\bfnm{B.}\binits{B.}},
  \bauthor{\bsnm{Hogg},~\bfnm{J.~S.}\binits{J.~S.}},
  \bauthor{\bsnm{Boissy},~\bfnm{R.}\binits{R.}},
  \bauthor{\bsnm{Yu},~\bfnm{S.}\binits{S.}},
  \bauthor{\bsnm{Powell},~\bfnm{E.}\binits{E.}},
  \bauthor{\bsnm{Keefe},~\bfnm{R.}\binits{R.}},
  \bauthor{\bsnm{Ehrlich},~\bfnm{N.~E.}\binits{N.~E.}},
  \bauthor{\bsnm{Shen},~\bfnm{K.}\binits{K.}},
  \bauthor{\bsnm{Hayes},~\bfnm{J.}\binits{J.}},
  \bauthor{\bsnm{Barbadora},~\bfnm{K.}\binits{K.}},
  \bauthor{\bsnm{Klimke},~\bfnm{W.}\binits{W.}},
  \bauthor{\bsnm{Dernovoy},~\bfnm{D.}\binits{D.}},
  \bauthor{\bsnm{Tatusova},~\bfnm{T.}\binits{T.}},
  \bauthor{\bsnm{Parkhill},~\bfnm{J.}\binits{J.}},
  \bauthor{\bsnm{Bentley},~\bfnm{S.~D.}\binits{S.~D.}},
  \bauthor{\bsnm{Post},~\bfnm{J.~C.}\binits{J.~C.}},
  \bauthor{\bsnm{Ehrlich},~\bfnm{G.~D.}\binits{G.~D.}} \AND
  \bauthor{\bsnm{Hu},~\bfnm{F.~Z.}\binits{F.~Z.}}
  (\byear{2007}).
\btitle{{C}omparative genomic analyses of seventeen
\textit{Streptococcus pneumoniae} strains: Insights into the pneumococcal
supragenome}.
\bjournal{J. Bacteriol.}
\bvolume{189}
\bpages{8186--8195}.
\end{barticle}
\endbibitem

\bibitem[\protect\citeauthoryear{Hogg et al.}{2007}]{pmid17550610}
\begin{barticle}[auto:SpringerTagBib|2009-01-14|16:51:27]
\bauthor{\bsnm{Hogg},~\bfnm{J.~S.}\binits{J.~S.}},
  \bauthor{\bsnm{Hu},~\bfnm{F.~Z.}\binits{F.~Z.}},
  \bauthor{\bsnm{Janto},~\bfnm{B.}\binits{B.}},
  \bauthor{\bsnm{Boissy},~\bfnm{R.}\binits{R.}},
  \bauthor{\bsnm{Hayes},~\bfnm{J.}\binits{J.}},
  \bauthor{\bsnm{Keefe},~\bfnm{R.}\binits{R.}},
  \bauthor{\bsnm{Post},~\bfnm{J.~C.}\binits{J.~C.}} \AND
  \bauthor{\bsnm{Ehrlich},~\bfnm{G.~D.}\binits{G.~D.}}
  (\byear{2007}).
\btitle{{C}haracterization and modeling of the
\textit{Haemophilus influenzae} core and supragenomes based on the complete
genomic sequences of {R}d and 12 clinical nontypeable strains}.
\bjournal{Genome Biol.}
\bvolume{8}
\bpages{R103}.
\end{barticle}
\endbibitem

\bibitem[\protect\citeauthoryear{Huson and
  Steel}{2004}]{Huson2004Bioinformatics15044248}
\begin{barticle}[pbm]
\bauthor{\bsnm{Huson},~\bfnm{Daniel~H.}\binits{D.~H.}} \AND
  \bauthor{\bsnm{Steel},~\bfnm{Mike}\binits{M.}}
(\byear{2004}).
\btitle{Phylogenetic trees based on gene content}.
\bjournal{Bioinformatics}
\bvolume{20}
\bpages{2044--2049}.
\bid{pmid={15044248}, doi={10.1093/bioinformatics/bth198}, pii={bth198}}
\end{barticle}
\endbibitem

\bibitem[\protect\citeauthoryear{Kettler et al.}{2007}]{pmid18159947}
\begin{barticle}[auto:SpringerTagBib|2009-01-14|16:51:27]
\bauthor{\bsnm{Kettler},~\bfnm{G.~C.}\binits{G.~C.}},
  \bauthor{\bsnm{Martiny},~\bfnm{A.~C.}\binits{A.~C.}},
  \bauthor{\bsnm{Huang},~\bfnm{K.}\binits{K.}},
  \bauthor{\bsnm{Zucker},~\bfnm{J.}\binits{J.}},
  \bauthor{\bsnm{Coleman},~\bfnm{M.~L.}\binits{M.~L.}},
  \bauthor{\bsnm{Rodrigue},~\bfnm{S.}\binits{S.}},
  \bauthor{\bsnm{Chen},~\bfnm{F.}\binits{F.}},
  \bauthor{\bsnm{Lapidus},~\bfnm{A.}\binits{A.}},
  \bauthor{\bsnm{Ferriera},~\bfnm{S.}\binits{S.}},
  \bauthor{\bsnm{Johnson},~\bfnm{J.}\binits{J.}},
  \bauthor{\bsnm{Steglich},~\bfnm{C.}\binits{C.}},
  \bauthor{\bsnm{Church},~\bfnm{G.~M.}\binits{G.~M.}},
  \bauthor{\bsnm{Richardson},~\bfnm{P.}\binits{P.}} \AND
  \bauthor{\bsnm{Chisholm},~\bfnm{S.~W.}\binits{S.~W.}}
  (\byear{2007}).
\btitle{Patterns and implications of gene gain and loss in
  the evolution of \textit{Prochlorococcus}}.
\bjournal{PLoS Genet.}
\bvolume{3}
\bpages{e231}.
\end{barticle}
\endbibitem

\bibitem[\protect\citeauthoryear{Kimura}{1964}]{Kimura1964}
\begin{barticle}[mr]
\bauthor{\bsnm{Kimura},~\bfnm{Motoo}\binits{M.}}
(\byear{1964}).
\btitle{Diffusion models in population genetics}.
\bjournal{J. Appl. Probab.}
\bvolume{1}
\bpages{177--232}.
\bid{mr={0172727}}
\end{barticle}
\endbibitem

\bibitem[\protect\citeauthoryear{Kingman}{1982}]{Kingman1982}
\begin{barticle}[mr]
\bauthor{\bsnm{Kingman},~\bfnm{J.~F.~C.}\binits{J.~F.~C.}}
(\byear{1982}).
\btitle{The coalescent}.
\bjournal{Stochastic Process. Appl.}
\bvolume{13}
\bpages{235--248}.
\bid{doi={10.1016/0304-4149(82)90011-4}, mr={671034}}
\end{barticle}
\endbibitem

\bibitem[\protect\citeauthoryear{Kunin and Ouzounis}{2003}]{pmid12874054}
\begin{barticle}[auto:SpringerTagBib|2009-01-14|16:51:27]
\bauthor{\bsnm{Kunin},~\bfnm{V.}\binits{V.}} \AND
  \bauthor{\bsnm{Ouzounis},~\bfnm{C.~A.}\binits{C.~A.}}
  (\byear{2003}).
\btitle{{G}ene{T}{R}{A}{C}{E}-reconstruction of gene content of
  ancestral species}.
\bjournal{Bioinformatics}
\bvolume{19}
\bpages{1412--1416}.
\end{barticle}
\endbibitem\vadjust{\goodbreak}

\bibitem[\protect\citeauthoryear{Lapierre and Gogarten}{2009}]{Lapierre2009}
\begin{barticle}[auto:SpringerTagBib|2009-01-14|16:51:27]
\bauthor{\bsnm{Lapierre},~\bfnm{P.}\binits{P.}} \AND
  \bauthor{\bsnm{Gogarten},~\bfnm{J.~P.}\binits{J.~P.}}
  (\byear{2009}).
\btitle{Estimating the size of the bacterial pan-genome}.
\bjournal{Trends in Genetics}
\bvolume{25}
\bpages{107--110}.
\end{barticle}
\endbibitem

\bibitem[\protect\citeauthoryear{Lef\'{e}bure and Stanhope}{2007}]{pmid17475002}
\begin{barticle}[auto:SpringerTagBib|2009-01-14|16:51:27]
\bauthor{\bsnm{Lef\'{e}bure},~\bfnm{T.}\binits{T.}} \AND
  \bauthor{\bsnm{Stanhope},~\bfnm{M.~J.}\binits{M.~J.}}
  (\byear{2007}).
\btitle{{E}volution of the core and pan-genome of
\textit{Streptococcus}: Positive selection, recombination, and genome
composition}.
\bjournal{Genome Biol.}
\bvolume{8}
\bpages{R71}.
\end{barticle}
\endbibitem

\bibitem[\protect\citeauthoryear{Maiden et al.}{1998}]{MaidenEtAl1998}
\begin{barticle}[auto:SpringerTagBib|2009-01-14|16:51:27]
\bauthor{\bsnm{Maiden},~\bfnm{M.~C.}\binits{M.~C.}},
  \bauthor{\bsnm{Bygraves},~\bfnm{J.~A.}\binits{J.~A.}},
  \bauthor{\bsnm{Feil},~\bfnm{E.}\binits{E.}},
  \bauthor{\bsnm{Morelli},~\bfnm{G.}\binits{G.}},
  \bauthor{\bsnm{Russell},~\bfnm{J.~E.}\binits{J.~E.}},
  \bauthor{\bsnm{Urwin},~\bfnm{R.}\binits{R.}},
  \bauthor{\bsnm{Zhang},~\bfnm{Q.}\binits{Q.}},
  \bauthor{\bsnm{Zhou},~\bfnm{J.}\binits{J.}},
  \bauthor{\bsnm{Zurth},~\bfnm{K.}\binits{K.}},
  \bauthor{\bsnm{Caugant},~\bfnm{D.~A.}\binits{D.~A.}},
  \bauthor{\bsnm{Feavers},~\bfnm{I.~M.}\binits{I.~M.}},
  \bauthor{\bsnm{Achtman},~\bfnm{M.}\binits{M.}} \AND
  \bauthor{\bsnm{Spratt},~\bfnm{B.~G.}\binits{B.~G.}}
  (\byear{1998}).
\btitle{Multilocus sequence typing: A portable approach to
the identification of clones within populations of pathogenic
microorganisms}.
\bjournal{Proc. Natl. Acad. Sci. USA}
\bvolume{95}
\bpages{3140--3145}.
\end{barticle}
\endbibitem

\bibitem[\protect\citeauthoryear{Maynard-Smith}{1995}]{MaynardSmith1995}
\begin{bincollection}[vtex]
\bauthor{\bsnm{Maynard-Smith},~\bfnm{J.}\binits{J.}}
(\byear{1995}).
\btitle{Do bacteria have population genetics}?
In \bbooktitle{Population Genetics of Bacteria}
\bpages{1--12}.
\bpublisher{Cambridge Univ. Press}, \baddress{Cambridge}.
\end{bincollection}
\endbibitem

\bibitem[\protect\citeauthoryear{Medini et al.}{2005}]{pmid16185861}
\begin{barticle}[pbm]
\bauthor{\bsnm{Medini},~\bfnm{Duccio}\binits{D.}},
  \bauthor{\bsnm{Donati},~\bfnm{Claudio}\binits{C.}},
  \bauthor{\bsnm{Tettelin},~\bfnm{Herv\'e}\binits{H.}},
  \bauthor{\bsnm{Masignani},~\bfnm{Vega}\binits{V.}} \AND
  \bauthor{\bsnm{Rappuoli},~\bfnm{Rino}\binits{R.}}
(\byear{2005}).
\btitle{The microbial pan-genome}.
\bjournal{Curr. Opin. Genet. Dev.}
\bvolume{15}
\bpages{589--594}.
\bid{pii={S0959-437X(05)00175-9}, doi={10.1016/j.gde.2005.09.006},
  pmid={16185861}}
\end{barticle}
\endbibitem

\bibitem[\protect\citeauthoryear{M{\"o}hle and
  Sagitov}{2001}]{MoehleSagitov2001}
\begin{barticle}[mr]
\bauthor{\bsnm{M{\"o}hle},~\bfnm{Martin}\binits{M.}} \AND
  \bauthor{\bsnm{Sagitov},~\bfnm{Serik}\binits{S.}}
(\byear{2001}).
\btitle{A classification of coalescent processes for haploid exchangeable
  population models}.
\bjournal{Ann. Probab.}
\bvolume{29}
\bpages{1547--1562}.
\bid{doi={10.1214/aop/1015345761}, mr={1880231}}
\end{barticle}
\endbibitem

\bibitem[\protect\citeauthoryear{Perna et al.}{2001}]{Perna2001}
\begin{barticle}[auto:SpringerTagBib|2009-01-14|16:51:27]
\bauthor{\bsnm{Perna},~\bfnm{N.~T.}\binits{N.~T.}},
  \bauthor{\bsnm{Plunkett},~\bfnm{G.}\binits{G.}},
  \bauthor{\bsnm{Burland},~\bfnm{V.}\binits{V.}},
  \bauthor{\bsnm{Mau},~\bfnm{B.}\binits{B.}},
  \bauthor{\bsnm{Glasner},~\bfnm{J.~D.}\binits{J.~D.}},
  \bauthor{\bsnm{Rose},~\bfnm{D.~J.}\binits{D.~J.}},
  \bauthor{\bsnm{Mayhew},~\bfnm{G.~F.}\binits{G.~F.}},
  \bauthor{\bsnm{Evans},~\bfnm{P.~S.}\binits{P.~S.}},
  \bauthor{\bsnm{Gregor},~\bfnm{J.}\binits{J.}},
  \bauthor{\bsnm{Kirkpatrick},~\bfnm{H.~A.}\binits{H.~A.}},
  \bauthor{\bsnm{P\'{e}sfai},~\bfnm{G.}\binits{G.}},
  \bauthor{\bsnm{Hackett},~\bfnm{J.}\binits{J.}},
  \bauthor{\bsnm{Klink},~\bfnm{S.}\binits{S.}},
  \bauthor{\bsnm{Boutin},~\bfnm{A.}\binits{A.}},
  \bauthor{\bsnm{Shao},~\bfnm{Y.}\binits{Y.}},
  \bauthor{\bsnm{Miller},~\bfnm{L.}\binits{L.}},
  \bauthor{\bsnm{Grotbeck},~\bfnm{E.~J.}\binits{E.~J.}},
  \bauthor{\bsnm{Davis},~\bfnm{N.~W.}\binits{N.~W.}},
  \bauthor{\bsnm{Lim},~\bfnm{A.}\binits{A.}},
  \bauthor{\bsnm{Dimalanta},~\bfnm{E.~T.}\binits{E.~T.}},
  \bauthor{\bsnm{Potamousis},~\bfnm{K.~D.}\binits{K.~D.}},
  \bauthor{\bsnm{Apodaca},~\bfnm{J.}\binits{J.}},
  \bauthor{\bsnm{Anantharaman},~\bfnm{T.~S.}\binits{T.~S.}},
  \bauthor{\bsnm{Lin},~\bfnm{J.}\binits{J.}},
  \bauthor{\bsnm{Yen},~\bfnm{G.}\binits{G.}},
  \bauthor{\bsnm{Schwartz},~\bfnm{D.~C.}\binits{D.~C.}},
  \bauthor{\bsnm{Welch},~\bfnm{R.~A.}\binits{R.~A.}} \AND
  \bauthor{\bsnm{Blattner},~\bfnm{F.~R.}\binits{F.~R.}}
  (\byear{2001}).
\btitle{Genome sequence of enterohaemorrhagic
\textit{Escherichia coli} {O}157:{H}7}.
\bjournal{Nature}
\bvolume{409}
\bpages{529--533}.
\end{barticle}
\endbibitem

\bibitem[\protect\citeauthoryear{Riley and Lizotte-Waniewski}{2009}]{Riley2009}
\begin{barticle}[auto:SpringerTagBib|2009-01-14|16:51:27]
\bauthor{\bsnm{Riley},~\bfnm{M.~A.}\binits{M.~A.}} \AND
  \bauthor{\bsnm{Lizotte-Waniewski},~\bfnm{M.}\binits{M.}}
  (\byear{2009}).
\btitle{Population genomics and the bacterial species
concept}.
\bjournal{Methods Mol. Biol.}
\bvolume{532}
\bpages{367--377}.
\end{barticle}
\endbibitem

\bibitem[\protect\citeauthoryear{Tettelin et al.}{2005}]{pmid16172379}
\begin{barticle}[auto:SpringerTagBib|2009-01-14|16:51:27]
\bauthor{\bsnm{Tettelin},~\bfnm{H.}\binits{H.}},
  \bauthor{\bsnm{Masignani},~\bfnm{V.}\binits{V.}},
  \bauthor{\bsnm{Cieslewicz},~\bfnm{M.~J.}\binits{M.~J.}},
  \bauthor{\bsnm{Donati},~\bfnm{C.}\binits{C.}},
  \bauthor{\bsnm{Medini},~\bfnm{D.}\binits{D.}},
  \bauthor{\bsnm{Ward},~\bfnm{N.~L.}\binits{N.~L.}},
  \bauthor{\bsnm{Angiuoli},~\bfnm{S.~V.}\binits{S.~V.}},
  \bauthor{\bsnm{Crabtree},~\bfnm{J.}\binits{J.}},
  \bauthor{\bsnm{Jones},~\bfnm{A.~L.}\binits{A.~L.}},
  \bauthor{\bsnm{Durkin},~\bfnm{A.~S.}\binits{A.~S.}},
  \bauthor{\bsnm{Deboy},~\bfnm{R.~T.}\binits{R.~T.}},
  \bauthor{\bsnm{Davidsen},~\bfnm{T.~M.}\binits{T.~M.}},
  \bauthor{\bsnm{Mora},~\bfnm{M.}\binits{M.}},
  \bauthor{\bsnm{Scarselli},~\bfnm{M.}\binits{M.}},
  \bauthor{\bparticle{Margarit~y }\bsnm{Ros},~\bfnm{I.}\binits{I.}},
  \bauthor{\bsnm{Peterson},~\bfnm{J.~D.}\binits{J.~D.}},
  \bauthor{\bsnm{Hauser},~\bfnm{C.~R.}\binits{C.~R.}},
  \bauthor{\bsnm{Sundaram},~\bfnm{J.~P.}\binits{J.~P.}},
  \bauthor{\bsnm{Nelson},~\bfnm{W.~C.}\binits{W.~C.}},
  \bauthor{\bsnm{Madupu},~\bfnm{R.}\binits{R.}},
  \bauthor{\bsnm{Brinkac},~\bfnm{L.~M.}\binits{L.~M.}},
  \bauthor{\bsnm{Dodson},~\bfnm{R.~J.}\binits{R.~J.}},
  \bauthor{\bsnm{Rosovitz},~\bfnm{M.~J.}\binits{M.~J.}},
  \bauthor{\bsnm{Sullivan},~\bfnm{S.~A.}\binits{S.~A.}},
  \bauthor{\bsnm{Daugherty},~\bfnm{S.~C.}\binits{S.~C.}},
  \bauthor{\bsnm{Haft},~\bfnm{D.~H.}\binits{D.~H.}},
  \bauthor{\bsnm{Selengut},~\bfnm{J.}\binits{J.}},
  \bauthor{\bsnm{Gwinn},~\bfnm{M.~L.}\binits{M.~L.}},
  \bauthor{\bsnm{Zhou},~\bfnm{L.}\binits{L.}},
  \bauthor{\bsnm{Zafar},~\bfnm{N.}\binits{N.}},
  \bauthor{\bsnm{Khouri},~\bfnm{H.}\binits{H.}},
  \bauthor{\bsnm{Radune},~\bfnm{D.}\binits{D.}},
  \bauthor{\bsnm{Dimitrov},~\bfnm{G.}\binits{G.}},
  \bauthor{\bsnm{Watkins},~\bfnm{K.}\binits{K.}},
  \bauthor{\bsnm{O'Connor},~\bfnm{K.~J.}\binits{K.~J.}},
  \bauthor{\bsnm{Smith},~\bfnm{S.}\binits{S.}},
  \bauthor{\bsnm{Utterback},~\bfnm{T.~R.}\binits{T.~R.}},
  \bauthor{\bsnm{White},~\bfnm{O.}\binits{O.}},
  \bauthor{\bsnm{Rubens},~\bfnm{C.~E.}\binits{C.~E.}},
  \bauthor{\bsnm{Grandi},~\bfnm{G.}\binits{G.}},
  \bauthor{\bsnm{Madoff},~\bfnm{L.~C.}\binits{L.~C.}},
  \bauthor{\bsnm{Kasper},~\bfnm{D.~L.}\binits{D.~L.}},
  \bauthor{\bsnm{Telford},~\bfnm{J.~L.}\binits{J.~L.}},
  \bauthor{\bsnm{Wessels},~\bfnm{M.~R.}\binits{M.~R.}},
  \bauthor{\bsnm{Rappuoli},~\bfnm{R.}\binits{R.}} \AND
  \bauthor{\bsnm{Fraser},~\bfnm{C.~M.}\binits{C.~M.}}
  (\byear{2005}).
\btitle{Genome analysis of multiple pathogenic isolates of
\textit{Streptococcus agalactiae}: Implications for the microbial
``pan-genome.''}
\bjournal{Proc. Natl. Acad. Sci. USA}
\bvolume{102}
\bpages{13950--13955}.
\end{barticle}
\endbibitem

\bibitem[\protect\citeauthoryear{Tettelin et al.}{2008}]{pmid19086349}
\begin{barticle}[pbm]
\bauthor{\bsnm{Tettelin},~\bfnm{Herv\'e}\binits{H.}},
  \bauthor{\bsnm{Riley},~\bfnm{David}\binits{D.}},
  \bauthor{\bsnm{Cattuto},~\bfnm{Ciro}\binits{C.}} \AND
  \bauthor{\bsnm{Medini},~\bfnm{Duccio}\binits{D.}}
(\byear{2008}).
\btitle{Comparative genomics: The bacterial pan-genome}.
\bjournal{Curr. Opin. Microbiol.}
\bvolume{11}
\bpages{472--477}.
\bid{pmid={19086349}}
\end{barticle}
\endbibitem

\bibitem[\protect\citeauthoryear{Vulic et al.}{1997}]{pmid9275198}
\begin{barticle}[auto:SpringerTagBib|2009-01-14|16:51:27]
\bauthor{\bsnm{Vulic},~\bfnm{M.}\binits{M.}},
\bauthor{\bsnm{Dionisio},~\bfnm{F.}\binits{F.}},
\bauthor{\bsnm{Taddei},~\bfnm{F.}\binits{F.}} \AND
\bauthor{\bsnm{Radman},~\bfnm{M.}\binits{M.}}
(\byear{1997}).
\btitle{Molecular keys to speciation: {D}{N}{A} polymorphism and the
  control of genetic exchange in enterobacteria}.
\bjournal{Proc. Natl. Acad. Sci. USA}
\bvolume{94}
\bpages{9763--9767}.
\end{barticle}
\endbibitem

\bibitem[\protect\citeauthoryear{Wakeley}{2008}]{Wakeley2008}
\begin{bbook}[auto:SpringerTagBib|2009-01-14|16:51:27]
\bauthor{\bsnm{Wakeley},~\bfnm{J.}\binits{J.}}
(\byear{2008}).
\btitle{Coalescent Theory: An Introduction}.
\bpublisher{Roberts and Company}, \baddress{Colorado}.
\end{bbook}
\endbibitem

\bibitem[\protect\citeauthoryear{Wright}{1938}]{Wright1938}
\begin{barticle}[auto:SpringerTagBib|2009-01-14|16:51:27]
\bauthor{\bsnm{Wright},~\bfnm{S.}\binits{S.}}
(\byear{1938}).
\btitle{The distribution of gene frequencies under irreversible
  mutation}.
\bjournal{Proc. Natl. Acad. Sci. USA}
\bvolume{24}
\bpages{253--259}.
\end{barticle}
\endbibitem

\end{thebibliography}
\end{document}